\documentclass[sort&compress,3p]{elsarticle}

\usepackage{amssymb}
\usepackage{amsmath}
\usepackage{mathtools}
\usepackage{bm}
\usepackage{multirow}
\usepackage{amsthm}
\usepackage{caption}
\usepackage{subcaption}
\usepackage{comment}
\usepackage{array}
\usepackage{xcolor}
\usepackage{graphicx}
\usepackage{multirow}

\newcommand{\modified}[1]{#1}

\def\N{\mathbb{N}}

\def\R{\mathbb{R}}

\def\BBname{Bernstein--B\'{e}zier}
\def\ie{i.e.}
\def\eg{e.g.}

\newcommand{\bfm}[1]{\bm{#1}}
\newcommand{\mc}[1]{\mathcal{#1}}

\newcommand{\bracket}[1]{\left( #1 \right)}
\newcommand{\set}[1]{\left\{ #1 \right\}}
\newcommand{\tset}[1]{\{ #1 \}}
\newcommand{\defset}[2]{\left\{ #1: \ #2 \right\}}
\newcommand{\abs}[1]{\left| #1 \right|}

\newcommand{\convh}[1]{[ #1 ]}
\newcommand{\C}[1]{{C}^{#1}}
\newcommand{\D}{\mathbb D}
\newcommand{\poly}{\mathbb P}
\newcommand{\tri}{\triangle}
\newcommand{\CTR}[1]{{#1}_{\mathrm{CT}}}
\newcommand{\CTT}{\CTR{\tri}}

\newcommand{\splspace}[3]{\mathbb{S}_{#1}^{#2}(#3)}

\newcommand{\grad}[1]{\nabla #1}

\newtheorem{remark}{Remark}

\newtheorem{example}{Example}

\begin{document}

\begin{frontmatter}


\address[FMF]{Faculty of Mathematics and Physics, University of Ljubljana, Jadranska 19, 1000 Ljubljana, Slovenia}
\address[IMFM]{Institute of Mathematics, Physics and Mechanics, Jadranska 19, 1000 Ljubljana, Slovenia}
\address[IJS]{Parallel and Distributed Systems Laboratory, ’Jo\v{z}ef Stefan’ Institute, Jamova cesta 39, 1000 Ljubljana, Slovenia}
\author[FMF]{Ema \v{C}e\v{s}ek}
\ead{ema.cesek@fmf.uni-lj.si}
\author[FMF,IMFM]{\corref{cor}Jan Gro\v{s}elj}
\ead{jan.groselj@fmf.uni-lj.si}
\author[IJS,FMF]{Andrej Kolar-Požun}
\ead{andrej.pozun@ijs.si}
\author[FMF,IMFM]{Maru\v{s}a Lek\v{s}e}
\ead{marusa.lekse@imfm.si}
\author[FMF]{Ga\v{s}per Domen Romih}
\ead{gasperdomen.romih@fmf.uni-lj.si}
\author[FMF,IMFM]{Ada \v{S}adl Praprotnik}
\ead{ada.sadl-praprotnik@imfm.si}
\author[FMF]{Matija \v{S}teblaj}
\ead{matija.steblaj@fmf.uni-lj.si}
\cortext[cor]{Corresponding author}

\title{A representation and comparison of three cubic macro-elements}


\begin{abstract}
The paper is concerned with three types of cubic splines over a triangulation that are characterized by three degrees of freedom associated with each vertex of the triangulation. The splines differ in computational complexity, polynomial reproduction properties, and smoothness. With the aim to make them a versatile tool for numerical analysis, a unified representation in terms of locally supported basis functions is established. The construction of these functions is based on geometric concepts and is expressed in the Bernstein--B\'{e}zier form. They are readily applicable in a range of standard approximation methods, which is demonstrated by a number of numerical experiments.
\end{abstract}

\begin{keyword}
splines over triangulations \sep Bernstein--B\'{e}zier techniques \sep macro-elements \sep spline basis \sep spline approximation
\end{keyword}

\end{frontmatter}


\section{Introduction}
\label{sec:intro}

Splines over triangulations are a widely used representation framework in numerical approximation methods. Thanks to the flexibility of triangulations as a meshing technique, they can be used to approximate data on complicated domains and, within a domain, by different levels of coarseness. In the classical sense, a spline over a triangulation is a function whose restriction to each triangle of the triangulation is a bivariate polynomial of a fixed total degree. At vertices and across edges of the triangles, these polynomials are glued together by a certain order of smoothness. High order of smoothness is desirable but not easy to achieve as smoothness constraints depend on geometry of the triangulation. For low polynomial degrees (relative to smoothness order) this prevents local and stable construction of general splines and affects their approximation power (see, \eg{}, \cite{lai_07}).

In this paper we consider three macro-elements, \ie{}, splines over a triangulation that are on each triangle uniquely determined by a limited number of degrees of freedom related to the region in the vicinity of the triangle. Traditionally, such splines are characterized by interpolation functionals associated with vertices of the triangulation and often (although undesirable) with additional points on the edges or inside the triangles (see, \eg{}, \cite{ciarlet_02,lai_07}). We focus on splines that are locally polynomials of total degree at most three, \ie{}, cubic splines, which are practically attractive on account of low evaluation cost. However, constructing a fast convergent and smooth approximation of this type is not trivial.

In order to construct a spline that is $\C{1}$-smooth at the vertices, it is natural to fix three degrees of freedom for each vertex which intrinsically represent the value and the gradient of the spline at the vertex. On a triangle of the triangulation this amounts to nine degrees of freedom if the three vertices of the triangle are considered. Since a bivariate cubic polynomial is determined by ten degrees of freedom, there is only one free parameter left for each triangle. This is not sufficient to impose global $\C{1}$-smoothness. Keeping the construction limited to the nine degrees of freedom corresponding to a triangle, it is not even possible to achieve cubic precision. On a more positive note, it is sufficient to ensure reproduction of quadratic polynomials, \eg{}, by using the classical Zienkiewicz element \cite{bazaley_65}, which is the first macro-element that we take into consideration.

If the spline construction does not ensure the reproduction of polynomials to the highest possible degree, we can not account for optimal convergence rates in approximation methods. To overcome this, we can sacrifice the complete locality of the construction on a triangle and use also the degrees of freedom corresponding to the neighboring triangles. Using the technique for imposing $\C{1}$-smoothness across an edge proposed in \cite{foley_92}, we develop a $\C{0}$-smooth macro-element that is, for a triangle with no edges on the boundary of the domain, determined by eighteen degrees of freedom.

Global $\C{1}$-smoothness of a macro-element is in approximation methods not always necessary and is indeed often renounced to keep the construction computationally simple and effective. However, it often has a positive effect on the accuracy of solutions and also provides visually more pleasant results. A common approach for achieving $\C{1}$-smoothness without increasing the spline degree is to refine the triangulation in such a way that each triangle is split into three smaller triangles. A classical Clough--Tocher spline \cite{ct_clough_65} over such a refinement has three degrees of freedom per vertex and additional degrees of freedom related to directional derivatives along the edges. It is possible to reduce these extra degrees of freedom by restricting the directional derivatives to be linear but such a constraint only allows reproduction of quadratic polynomials. We instead follow the construction introduced in \cite{ct_mann_99} that ensures full reproduction. This provides a $\C{1}$-smooth macro-element with cubic precision that is based on the same amount of data as the second macro-element.

With the aim to use the described constructions in approximation methods, we establish a unified representation of the three macro-elements in terms of locally supported basis functions. They are defined and analyzed using the Bernstein--B\'{e}zier techniques. For each vertex, there are three basis functions that are uniquely determined by a choice of a triangle associated with the vertex. This geometric approach reflects recent ideas in constructing B-spline-like bases for other common macro-elements (see, \eg{}, \cite{ps_dierckx_97,ct3_speleers_10,ps12_lamnii_15,ps3_groselj_17,ct_groselj_22}).

Finally, we provide a numerical comparison of the macro-elements. The proposed representation is utilized in function approximation and solving boundary value problems. We perform tests of the best $L^2$ function approximation and scattered data approximation. In the context of boundary value problems, we consider the model Poisson problem and conduct examples showing the behavior of the macro-elements in the isoparametric finite element method and immersed penalized boundary method.

The remainder of the paper is organized as follows. In Section~\ref{sec:bbrep} we review relevant Bernstein--B\'{e}zier techniques for representing cubic splines and imposing $\C{1}$-smoothness. In Section~\ref{sec:cubic_splines} we provide the construction of the three macro-elements and define basis functions associated with the vertices of the triangulation. Section~\ref{sec:examples} contains numerical examples. We conclude the paper with some remarks in Section~\ref{sec:conclusion}.

\section{Bernstein--B\'{e}zier spline techniques}
\label{sec:bbrep}

The aim of this section is to present some preliminaries related to the construction of cubic splines over a triangulation in terms of the \BBname{} form. We give a particular emphasis on the constraints for (local) $\C{1}$-smoothness that is tackled \modified{in} Section~\ref{sec:cubic_splines}.

\subsection{Cubic splines over a triangulation}

Let $\Theta \subset \R^2$ be a polygonal domain. A triangulation $\tri = (V, E, T)$ of $\Theta$ is a discretization of the domain determined by a finite set of vertices $V$, a finite set of edges $E$, and a finite set of triangles $T$.

More precisely, for $n_v \in \N$, the set $V = \tset{\bfm{v}_1, \bfm{v}_2, \ldots, \bfm{v}_{n_v}}$ is a set of points in $\Theta$. The set $E = \set{\bfm{e}_1, \bfm{e}_2, \ldots, \bfm{e}_{n_e}}$ consists of $n_e \in \N$ line segments, each of the form $\bfm{e}_l = \convh{\bfm{v}_i, \bfm{v}_j}$ for a pair of two distinct vertices $\bfm{v}_i, \bfm{v}_j \in V$.
% The intersection of two distinct edges from $E$ is either empty or a point in $V$.
Similarly, the set $T = \set{\bfm{t}_1, \bfm{t}_2, \ldots, \bfm{t}_{n_t}}$ consists of $n_t \in \N$ triangles, each of the form $\bfm{t}_m = \convh{\bfm{v}_i, \bfm{v}_j, \bfm{v}_k}$ for a triplet of three affinely independent vertices $\bfm{v}_i, \bfm{v}_j, \bfm{v}_k \in V$.
% The edges of a triangle from $T$ are in $E$,
Vice versa, $V$ and $E$ are the set of vertices and the set of edges of all triangles in $T$, respectively.
We impose that the intersection of two distinct triangles from $T$ is empty, a point in $V$, or a line segment in $E$, and the union of triangles in $T$ is $\Theta$. Additionally, we denote by $E^b \subseteq E$ the set of cardinality $n_e^b \in \N$ comprising all edges that lie on the boundary of $\Theta$. 

The vector space of cubic splines over a triangulation $\tri$ of smoothness order $r \in \set{0,1}$ is defined by
\begin{equation*}
\splspace{3}{r}{\tri} \modified{\coloneqq} \defset{S \in C^r(\Theta)}{S|_{\bfm{t}_m} \in \poly_3, \ \bfm{t}_m \in T},
\end{equation*}
where\modified{, in general, $\poly_d$} denotes the space of bivariate polynomials of total degree at most \modified{$d \in \N_0$}. The dimension of $\splspace{3}{0}{\tri}$ is $n_v + 2 n_e + n_t$. Expressing the dimension of $\splspace{3}{1}{\tri}$ is much more complicated and is theoretically not fully resolved. It is known that
\begin{equation}
\label{eq:c1_dim_lb}
\dim(\splspace{3}{1}{\tri}) \geq 3 n_v^b + 2 (n_v - n_v^b) +1,
\end{equation}
where $n_v^b \in \N$ is the number of points in $V$ that lie on the boundary of $\Theta$. Moreover, the dimension depends on the geometry of the triangulation and is larger than the lower bound in \eqref{eq:c1_dim_lb} if the triangulation has vertices in the interior of $\Theta$ in which the edges meet with only two different slopes (see, \eg{}, \cite{lai_07}, for more details). Therefore, it is unclear how to construct and use such splines in approximation methods. In Section~\ref{sec:cubic_splines} we consider related spline spaces that can be more easily characterized.

\subsection{Bernstein--B\'{e}zier form}

As the restriction $S|_{\bfm{t}_m}$ of a spline $S \in \splspace{3}{0}{\tri}$ to a triangle $\bfm{t}_m = \convh{\bfm{v}_i, \bfm{v}_j, \bfm{v}_k} \in T$ is an element of $\poly_3$, it can be uniquely expressed in the cubic Bernstein basis. Let
\begin{equation*}
\D_3 \modified{\coloneqq} \defset{(\delta_0, \delta_1, \delta_2) \in \N_0^3}{\delta_0 + \delta_1 + \delta_2 = 3}
\end{equation*}
be a set consisting of ten multi-indices. For $\bfm{d} = (d_0, d_1, d_2) \in \D_3$ we define the Bernstein basis polynomial $B_{m,\bfm{d}} \in \poly_3$ associated with the triangle $\bfm{t}_m$ as
\begin{equation*}
B_{m,\bfm{d}}(\bfm{p}) \modified{\coloneqq} \frac{3!}{d_0! \, d_1! \, d_2!} \bracket{\tau_{m,i}(\bfm{p})}^{d_0} \bracket{\tau_{m,j}(\bfm{p})}^{d_1} \bracket{\tau_{m,k}(\bfm{p})}^{d_2},
\end{equation*}
where $(\tau_{m,i}(\bfm{p}), \tau_{m,j}(\bfm{p}), \tau_{m,k}(\bfm{p})) \in \R^3$ denotes the barycentric coordinates of a point $\bfm{p} \in \R^2$ with respect to $\bfm{t}_m$, i.e., \modified{the} unique triplet satisfying
\begin{equation*}
1 = \tau_{m,i}(\bfm{p}) + \tau_{m,j}(\bfm{p}) + \tau_{m,k}(\bfm{p}), \quad
\bfm{p} = \tau_{m,i}(\bfm{p}) \bfm{v}_i + \tau_{m,j}(\bfm{p}) \bfm{v}_j + \tau_{m,k}(\bfm{p}) \bfm{v}_k.
\end{equation*}
It is well known that the polynomial $S|_{\bfm{t}_m}$ can be uniquely represented in the Bernstein--B\'{e}zier form
\begin{equation*}
S|_{\bfm{t}_m} = \sum_{\bfm{d} \in \D_3} \beta_{m,\bfm{d}}(S) \, B_{m,\bfm{d}},
\end{equation*}
where $\beta_{m,\bfm{d}}: \modified{\splspace{3}{0}{\tri}} \rightarrow \R$ \modified{applied to $S$ corresponds to} the dual functional of $B_{m, \bfm{d}}$ \modified{applied to $S|_{\bfm{t}_m}$}. We refer the reader to, \eg{}, \cite{lai_07} for more details on this topic.

\subsection{Techniques for imposing $\C{1}$-smoothness}
\label{sec:c1_smoothness}

A spline $S \in \splspace{3}{0}{\tri}$ is $\C{1}$-smooth at a vertex $\bfm{v}_i \in V$ if and only if there exists a polynomial $P \in \poly_3$ such that $S(\bfm{v}_i) = P(\bfm{v}_i)$ and $\grad{S|_{\bfm{t}_m}}(\bfm{v}_i) = \grad{P}(\bfm{v}_i)$ for every triangle $\bfm{t}_m \in T$ with a vertex at $\bfm{v}_i$. Here and hereafter, $\grad$ denotes the gradient operator. Let us consider the first order Taylor polynomial
\begin{equation*}
Q_i(\bfm{p}) = P(\bfm{v}_i) + \grad{P}(\bfm{v}_i) \bfm{\cdot} (\bfm{p} - \bfm{v}_i),
\end{equation*}
where $\bfm{\cdot}$ denotes the inner product of two vectors. As $Q_i(\bfm{v}_i) = P(\bfm{v}_i)$ and $\grad{Q_i}(\bfm{v}_i) = \grad{P}(\bfm{v}_i)$, the standard vertex interpolation formulas for the \BBname{} form (see, \eg{}, \cite{lai_07}) imply
\begin{equation}
\label{eq:vertex_interpol}
\beta_{m,(3,0,0)}(S) = Q_i(\bfm{v}_i), \quad
\beta_{m,(2,1,0)}(S) = Q_i(\tfrac{2}{3} \bfm{v}_i + \tfrac{1}{3} \bfm{v}_j), \quad
\beta_{m,(2,0,1)}(S) = Q_i(\tfrac{2}{3} \bfm{v}_i + \tfrac{1}{3} \bfm{v}_k),
\end{equation}
where $\bfm{t}_m = \convh{\bfm{v}_i, \bfm{v}_j, \bfm{v}_k}$. This provides a convenient alternative characterization of $\C{1}$-smoothness at a vertex, \ie{}, $S$ is $\C{1}$-smooth at $\bfm{v}_i$ if and only if there exists a linear polynomial $Q_i$ such that \eqref{eq:vertex_interpol} holds for every triangle $\bfm{t}_m \in T$ with a vertex at $\bfm{v}_i$. Hence, by fixing a linear polynomial for each vertex, all but one coefficient in the Bernstein--B\'ezier form of $S|_{\bfm{t}_m}$ are uniquely determined.

\begin{comment}
A spline $S \in \splspace{3}{0}{\tri}$ is $\C{1}$-smooth at a vertex $\bfm{v}_i \in V$ if and only if there exists a linear bivariate polynomial $Q_i$ such that
\begin{equation}
\label{eq:vertex_interpol}
\beta_{m,(3,0,0)}(S) = Q_i(\bfm{v}_i), \quad
\beta_{m,(2,1,0)}(S) = Q_i(\tfrac{2}{3} \bfm{v}_i + \tfrac{1}{3} \bfm{v}_j), \quad
\beta_{m,(2,0,1)}(S) = Q_i(\tfrac{2}{3} \bfm{v}_i + \tfrac{1}{3} \bfm{v}_k)
\end{equation}
for every triangle $\bfm{t}_m = \convh{\bfm{v}_i, \bfm{v}_j, \bfm{v}_k} \in T$ with a vertex at $\bfm{v}_i$. Hence, by fixing a linear polynomial for each vertex, all but one coefficient in the Bernstein--B\'ezier form of $S|_{\bfm{t}_m}$ are uniquely determined.
\end{comment}

Let us consider triangles $\bfm{t}_m = \convh{\bfm{v}_i, \bfm{v}_j, \bfm{v}_k} \in T$ and $\bfm{t}_{m'} = \convh{\bfm{v}_i, \bfm{v}_j, \bfm{v}_{k'}} \in T$ that have the edge $\bfm{e}_l = \convh{\bfm{v}_i, \bfm{v}_j} \in E \setminus E^b$ in common. Provided that $S$ is $\C{1}$-smooth at the vertices $\bfm{v}_i$ and $\bfm{v}_j$, the $\C{1}$-smoothness across $\convh{\bfm{v}_i, \bfm{v}_j}$ can be achieved by selecting $\beta_{m,(1,1,1)}(S)$ and $\beta_{m',(1,1,1)}(S)$ such that the value of
\begin{equation}
\label{eq:c1_condition_edge}
\xi_l(S) \modified{\coloneqq}
\tau_{m,i}(\bfm{v}_{k'}) \beta_{m,(2,1,0)}(S) +
\tau_{m,j}(\bfm{v}_{k'}) \beta_{m,(1,2,0)}(S) +
\tau_{m,k}(\bfm{v}_{k'}) \beta_{m,(1,1,1)}(S) -
\beta_{m',(1,1,1)}(S)
\end{equation}
is equal to zero. However, imposing $\C{1}$-smoothness across all interior edges of $\tri$ is challenging as it results in a global system of equations depending on the geometry of the triangulation. In fact, since for certain triangulations the lower bound \eqref{eq:c1_dim_lb} is attained, fixing three degrees of freedom per vertex implies that such a global system in general does not have a solution.

A convenient technique for imposing $\C{1}$-smoothness across a single edge $\bfm{e}_l$ without selecting an additional parameter was proposed in \cite{foley_92}.
If $S|_{\bfm{t}_m \cup \bfm{t}_{m'}}$ is $\C{2}$-smooth, then
\begin{subequations}
\label{eq:foley_opitz_c2_conditions}
\begin{equation}
\begin{split}
\beta_{m',(1,0,2)}(S) &= \tau_{m,i}(\bfm{v}_{k'})^2 \beta_{m,(3,0,0)}(S) + 2 \tau_{m,i}(\bfm{v}_{k'}) \tau_{m,j}(\bfm{v}_{k'}) \beta_{m,(2,1,0)}(S) \\
&\quad + \tau_{m,j}(\bfm{v}_{k'})^2 \beta_{m,(1,2,0)}(S) + 2 \tau_{m,i}(\bfm{v}_{k'}) \tau_{m,k}(\bfm{v}_{k'}) \beta_{m,(2,0,1)}(S) \\
&\quad + 2 \tau_{m,j}(\bfm{v}_{k'}) \tau_{m,k}(\bfm{v}_{k'}) \beta_{m,(1,1,1)}(S) + \tau_{m,k}(\bfm{v}_{k'})^2 \beta_{m,(1,0,2)}(S)
\end{split}
\end{equation}
and
\begin{equation}
\begin{split}
\beta_{m',(0,1,2)}(S) &= \tau_{m,i}(\bfm{v}_{k'})^2 \beta_{m,(2,1,0)}(S) + 2 \tau_{m,i}(\bfm{v}_{k'}) \tau_{m,j}(\bfm{v}_{k'}) \beta_{m,(1,2,0)}(S) \\
&\quad + \tau_{m,j}(\bfm{v}_{k'})^2 \beta_{m,(0,3,0)}(S) + 2 \tau_{m,i}(\bfm{v}_{k'}) \tau_{m,k}(\bfm{v}_{k'}) \beta_{m,(1,1,1)}(S) \\
&\quad + 2 \tau_{m,j}(\bfm{v}_{k'}) \tau_{m,k}(\bfm{v}_{k'}) \beta_{m,(0,2,1)}(S) + \tau_{m,k}(\bfm{v}_{k'})^2 \beta_{m,(0,1,2)}(S).
\end{split}
\end{equation}
\end{subequations}
All the values appearing in \eqref{eq:foley_opitz_c2_conditions} are fixed by vertex interpolation, i.e., by expressions analogous to \eqref{eq:vertex_interpol}, except for $\beta_{m,(1,1,1)}(S)$. Motivated by this, we define the functional
\begin{equation}
\label{eq:gamma}
\begin{split}
\gamma_{m,k}(S) &\modified{\coloneqq}
\frac{1}{2 \tau_{m,k}(\bfm{v}_{k'}) (1 - \tau_{m,k}(\bfm{v}_{k'}))} \left(
\beta_{m',(1,0,2)}(S) + \beta_{m',(0,1,2)}(S) \right. \\[0.2em]
&\quad - \tau_{m,i}(\bfm{v}_{k'})^2 (\beta_{m,(3,0,0)}(S) + \beta_{m,(2,1,0)}(S)) \\
&\quad - 2 \tau_{m,i}(\bfm{v}_{k'}) \tau_{m,j}(\bfm{v}_{k'}) (\beta_{m,(2,1,0)}(S) + \beta_{m,(1,2,0)}(S)) \\[0.2em]
&\quad - \tau_{m,j}(\bfm{v}_{k'})^2 (\beta_{m,(1,2,0)}(S) + \beta_{m,(0,3,0)}(S)) \\[0.2em]
&\quad - 2 \tau_{m,i}(\bfm{v}_{k'}) \tau_{m,k}(\bfm{v}_{k'}) \beta_{m,(2,0,1)}(S) \\[0.2em]
&\quad - 2 \tau_{m,j}(\bfm{v}_{k'}) \tau_{m,k}(\bfm{v}_{k'}) \beta_{m,(0,2,1)}(S) \\[0.2em]
&\quad \left. - \tau_{m,k}(\bfm{v}_{k'})^2 (\beta_{m,(1,0,2)}(S) + \beta_{m,(0,1,2)}(S)) \right),
\end{split}
\end{equation}
which is obtained by combining the equations in \eqref{eq:foley_opitz_c2_conditions} and expressing $\beta_{m,(1,1,1)}(S)$.

Setting $\beta_{m,(1,1,1)}(S) = \gamma_{m,k}(S)$ and $\beta_{m',(1,1,1)}(S) = \gamma_{m',k'}(S)$ has some immediate benefits. Even if the assumptions leading to \eqref{eq:foley_opitz_c2_conditions} do not hold, such choice ensures that $\xi_l(S) = 0$, \ie{}, $S$ is $\C{1}$-smooth across $\bfm{e}_l$. Moreover, if the interpolation values at the vertices $\bfm{v}_i$, $\bfm{v}_j$, $\bfm{v}_k$, $\bfm{v}_{k'}$ come from the same polynomial $P \in \poly_3$, then $S|_{\bfm{t}_m \cup \bfm{t}_{m'}} = P$. For the proofs of these statements we refer the reader to \cite{foley_92}.

\section{Cubic macro-elements}
\label{sec:cubic_splines}

In this section we present three constructions of cubic splines over a triangulation $\triangle$\modified{, each of them determined by three degrees of freedom associated with every vertex of $\triangle$}. The splines are specified by providing the \BBname{} form on each triangle of $\triangle$ in terms of linear polynomials $Q_i$ associated with the vertices $\bfm{v}_i \in V$. \modified{For the first and third construction this is a reinterpretation of the results derived in \cite{bazaley_65} and \cite{ct_mann_99}, respectively, whereas the second construction appears to be new.} In the second part of the section we \modified{propose} locally supported basis functions for a numerically convenient representation of the considered splines.

\subsection{$\C{0}$-smooth macro-element with quadratic precision}
\label{sec:s1}

First, we recast the Zienkiewicz element \cite{bazaley_65} that provides a $C^0$-smooth spline construction over $\triangle$ with quadratic precision.

Let $\bfm{t}_m = [\bfm{v}_i, \bfm{v}_j, \bfm{v}_k] \in T$ be an arbitrary triangle of $\triangle$. We prescribe a spline $S \in \splspace{3}{0}{\triangle}$ by setting
\begin{subequations}
\label{eq:s1_vertex_interpolation}
\begin{align}
\beta_{m,(3,0,0)}(S) &= Q_i(\bfm{v}_i), &
\beta_{m,(2,1,0)}(S) &= Q_i(\tfrac{2}{3} \bfm{v}_i + \tfrac{1}{3} \bfm{v}_j), &
\beta_{m,(2,0,1)}(S) &= Q_i(\tfrac{2}{3} \bfm{v}_i + \tfrac{1}{3} \bfm{v}_k), \\
\beta_{m,(0,3,0)}(S) &= Q_j(\bfm{v}_j), &
\beta_{m,(0,2,1)}(S) &= Q_j(\tfrac{2}{3} \bfm{v}_j + \tfrac{1}{3} \bfm{v}_k), &
\beta_{m,(1,2,0)}(S) &= Q_j(\tfrac{2}{3} \bfm{v}_j + \tfrac{1}{3} \bfm{v}_i), \\
\beta_{m,(0,0,3)}(S) &= Q_k(\bfm{v}_k), &
\beta_{m,(1,0,2)}(S) &= Q_k(\tfrac{2}{3} \bfm{v}_k + \tfrac{1}{3} \bfm{v}_i), &
\beta_{m,(0,1,2)}(S) &= Q_k(\tfrac{2}{3} \bfm{v}_k + \tfrac{1}{3} \bfm{v}_j),
\end{align}
\end{subequations}
and
\begin{equation}
\label{eq:s1_central_value}
\beta_{m,(1,1,1)}(S) =
\tfrac{1}{3} Q_i(\tfrac{1}{2} \bfm{v}_i + \tfrac{1}{4} \bfm{v}_j + \tfrac{1}{4} \bfm{v}_k) +
\tfrac{1}{3} Q_j(\tfrac{1}{2} \bfm{v}_j + \tfrac{1}{4} \bfm{v}_k + \tfrac{1}{4} \bfm{v}_i) +
\tfrac{1}{3} Q_k(\tfrac{1}{2} \bfm{v}_k + \tfrac{1}{4} \bfm{v}_i + \tfrac{1}{4} \bfm{v}_j).
\end{equation}
Following the discussion in Section~\ref{sec:c1_smoothness}, such a spline is $\C{1}$-smooth at each vertex of $\triangle$ but only $\C{0}$-smooth across the edges of $\triangle$.

We denote by $\mathbb{S}_1 \subset \splspace{3}{0}{\triangle}$ the vector space consisting of all splines that can be obtained by the described construction. It holds that $\dim(\mathbb{S}_1) = 3 n_v$ and $\poly_2 \subset \mathbb{S}_1$. % where $\poly_2$ denotes the space of bivariate polynomials of total degree at most $2$.

\subsection{$\C{0}$-smooth macro-element with cubic precision}
\label{sec:s2}

Next, we describe the Bernstein--B\'{e}zier form on a triangle $\bfm{t}_m = [\bfm{v}_i, \bfm{v}_j, \bfm{v}_k] \in T$ that provides a $C^0$-smooth spline construction over $\triangle$ with cubic precision. For this we assume that $\bfm{t}_m$ has at most two edges in $E^b$. Additionally, we choose the ordering of the vertices of $\bfm{t}_m$ in such a way that $\convh{\bfm{v}_i, \bfm{v}_j} \in E \setminus E^b$ and, if $\bfm{t}_m$ has at most one boundary edge, $\convh{\bfm{v}_j, \bfm{v}_k} \in E \setminus E^b$.

We prescribe a spline $S \in \splspace{3}{0}{\triangle}$ by setting $\beta_{m,\bfm{d}}(S)$, $\bfm{d} \in \D_3 \setminus \set{(1,1,1)}$, as in \eqref{eq:s1_vertex_interpolation} and
\begin{equation*}
\beta_{m,(1,1,1)}(S) =
\begin{cases}
\frac{1}{3} \gamma_{m,k}(S) + \frac{1}{3} \gamma_{m,i}(S) + \frac{1}{3} \gamma_{m,j}(S) & \text{if all edges of $\bfm{t}_m$ are in $E \setminus E^b$,} \\
\frac{1}{2} \gamma_{m,k}(S) + \frac{1}{2} \gamma_{m,i}(S) & \text{if two edges of $\bfm{t}_m$ are in $E \setminus E^b$,} \\
\gamma_{m,k}(S) & \text{if one edge of $\bfm{t}_m$ is in $E \setminus E^b$,}
\end{cases}
\end{equation*}
where the functionals $\gamma_{m,i}$, $\gamma_{m,j}$, $\gamma_{m,k}$ are defined as in \eqref{eq:gamma}. Following the discussion in Section~\ref{sec:c1_smoothness}, the resulting spline is again $\C{1}$-smooth at each vertex of $\triangle$ and $\C{0}$-smooth across the edges of $\triangle$. The main difference in comparison with the spline given in Section~\ref{sec:s1} is that the value $\beta_{m,(1,1,1)}(S)$ depends not only on the polynomials $Q_i$, $Q_j$, $Q_k$ but also on the polynomials associated with the vertices of the triangles in $T$ that share an edge with $\bfm{t}_m$. Nonetheless, the construction remains local.

We denote by $\mathbb{S}_2 \subset \splspace{3}{0}{\triangle}$ the vector space consisting of all splines that can be obtained by the described construction. It holds that $\dim(\mathbb{S}_2) = 3 n_v$ and $\poly_3 \subset \mathbb{S}_2$.

\subsection{$\C{1}$-smooth macro-element with cubic precision}
\label{sec:s3}

\begin{figure}[!t]
\centering

\begin{subfigure}[b]{0.47\textwidth}
\centering
\includegraphics[width=\textwidth]{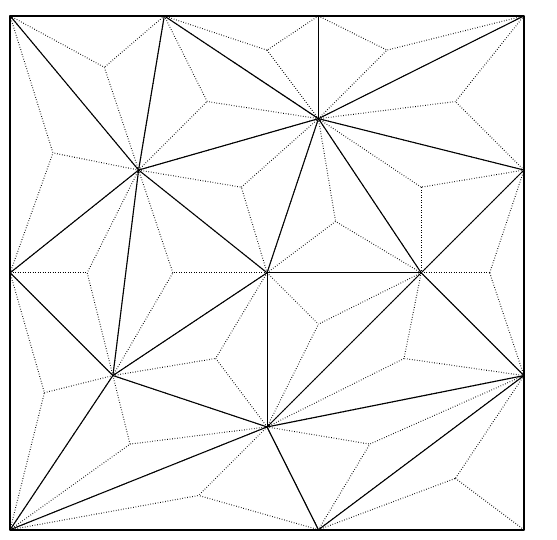}
\end{subfigure}
~
\begin{subfigure}[b]{0.47\textwidth}
\centering
\includegraphics[width=\textwidth]{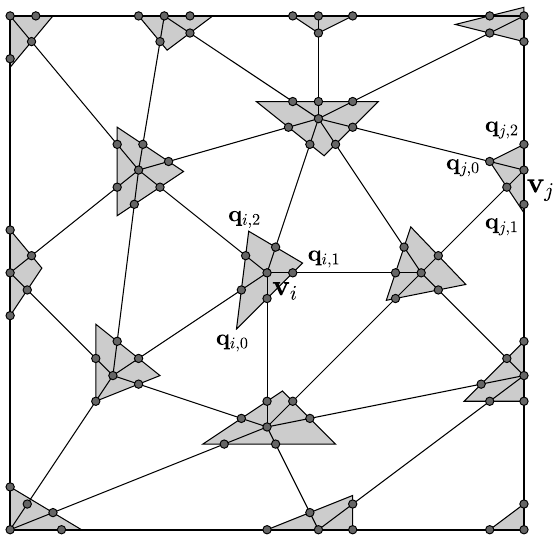}
\end{subfigure}

\caption{A triangulation $\tri$ of the unit square depicted by the uninterrupted lines. The figure on the left also shows the Clough--Tocher refinement $\CTT$ of $\tri$ determined by the dotted lines. The figure on the right depicts a configuration of triangles associated with the vertices of $\triangle$ suitable for construction of basis functions. Each triangle is obtained as the smallest triangle that contains the gray colored points in the vicinity of the vertex.}
\label{fig:triangulation}
\end{figure}

In order to obtain a $\C{1}$-smooth spline, a standard macro-element approach is to split each triangle of the triangulation into a number of smaller ones. Here, we use the Clough--Tocher splitting \cite{ct_clough_65} which refines a triangle $\bfm{t}_m = \convh{\bfm{v}_i, \bfm{v}_j, \bfm{v}_k} \in T$ into three triangles
\begin{equation*}
\bfm{t}_{m_k} = \convh{\bfm{v}_i, \bfm{v}_j, \bfm{v}_{ijk}}, \quad \bfm{t}_{m_i} = \convh{\bfm{v}_j, \bfm{v}_k, \bfm{v}_{ijk}}, \quad \bfm{t}_{m_j} = \convh{\bfm{v}_k, \bfm{v}_i, \bfm{v}_{ijk}}
\end{equation*}
determined by the split point $\bfm{v}_{ijk} = \frac{1}{3} \bfm{v}_i + \frac{1}{3} \bfm{v}_j + \frac{1}{3} \bfm{v}_k$. These triangles define a finer triangulation denoted by $\CTT$. For an illustration, see Figure~\ref{fig:triangulation} (left).

In what follows we provide a $\C{1}$-smooth spline construction over $\CTT$ \modified{inspired by \cite{ct_mann_99}}. It differs from the classical one in that it does not introduce extra degrees of freedom related to directional derivatives along the edges. As it was already proposed in \cite{ct_clough_65}, these degrees of freedom can be simply avoided by imposing that the directional derivatives are linear but such a constraint reduces precision of splines to quadratic polynomials. \modified{On the other hand,} the concepts introduced in \cite{ct_mann_99} \modified{ensure} cubic precision.

We prescribe a spline $S \in \splspace{3}{1}{\CTT}$ by first setting
\begin{align*}
\beta_{m_k,(3,0,0)}(S) &= Q_i(\bfm{v}_i), &
\beta_{m_k,(2,1,0)}(S) &= Q_i(\tfrac{2}{3} \bfm{v}_i + \tfrac{1}{3} \bfm{v}_j), &
\beta_{m_k,(2,0,1)}(S) &= Q_i(\tfrac{2}{3} \bfm{v}_i + \tfrac{1}{3} \bfm{v}_{ijk}), \\
\beta_{m_k,(0,3,0)}(S) &= Q_j(\bfm{v}_j), &
\beta_{m_k,(1,2,0)}(S) &= Q_j(\tfrac{2}{3} \bfm{v}_j + \tfrac{1}{3} \bfm{v}_i), &
\beta_{m_k,(0,2,1)}(S) &= Q_j(\tfrac{2}{3} \bfm{v}_j + \tfrac{1}{3} \bfm{v}_{ijk}),
\end{align*}
and, assuming $\bfm{t}_m$ has at most two edges in $E^b$,
\begin{align*}
\beta_{m_k,(1,1,1)}(S) &=
\frac{1}{3} Q_i(\tfrac{2}{3} \bfm{v}_i + \tfrac{1}{3} \bfm{v}_j)
+ \frac{1}{3} Q_j(\tfrac{2}{3} \bfm{v}_j + \tfrac{1}{3} \bfm{v}_i) \\
&\quad + \frac{1}{3}
\begin{cases}
\gamma_{m,k}(S) & \text{if $\convh{\bfm{v}_i, \bfm{v}_j} \in E \setminus E^b$,} \\
\frac{1}{2} \gamma_{m,i}(S) + \frac{1}{2} \gamma_{m,j}(S) & \text{if $\convh{\bfm{v}_i, \bfm{v}_j} \in E^b$ and $\convh{\bfm{v}_j, \bfm{v}_k}, \convh{\bfm{v}_k, \bfm{v}_i} \in E \setminus E^b$,} \\
\gamma_{m,j}(S) & \text{if $\convh{\bfm{v}_i, \bfm{v}_j}, \convh{\bfm{v}_j, \bfm{v}_k} \in E^b$ and $\convh{\bfm{v}_k, \bfm{v}_i} \in E \setminus E^b$,} \\
\gamma_{m,i}(S) & \text{if $\convh{\bfm{v}_i, \bfm{v}_j}, \convh{\bfm{v}_k, \bfm{v}_i} \in E^b$ and $\convh{\bfm{v}_j, \bfm{v}_k} \in E \setminus E^b$.}
\end{cases}
\end{align*}
We analogously define the values $\beta_{m_j,(d_0, d_1, d_2)}(S)$ and $\beta_{m_i,(d_0,d_1,d_2)}(S)$ for $d_2 \leq 1$. The remaining values are given by
\begin{align*}
\beta_{m_k,(1,0,2)}(S) = \beta_{m_j,(0,1,2)}(S) = \tfrac{1}{3} \beta_{m_k,(2,0,1)}(S) + \tfrac{1}{3} \beta_{m_k,(1,1,1)}(S) + \tfrac{1}{3} \beta_{m_j,(1,1,1)}(S), \\
\beta_{m_i,(1,0,2)}(S) = \beta_{m_k,(0,1,2)}(S) = \tfrac{1}{3} \beta_{m_i,(2,0,1)}(S) + \tfrac{1}{3} \beta_{m_i,(1,1,1)}(S) + \tfrac{1}{3} \beta_{m_k,(1,1,1)}(S), \\
\beta_{m_j,(1,0,2)}(S) = \beta_{m_i,(0,1,2)}(S) = \tfrac{1}{3} \beta_{m_j,(2,0,1)}(S) + \tfrac{1}{3} \beta_{m_j,(1,1,1)}(S) + \tfrac{1}{3} \beta_{m_i,(1,1,1)}(S),
\end{align*}
and
\begin{equation*}
\beta_{m_k,(0,0,3)}(S) = \beta_{m_i,(0,0,3)}(S) = \beta_{m_j,(0,0,3)}(S) = \tfrac{1}{3} \beta_{m_k,(1,0,2)}(S) + \tfrac{1}{3} \beta_{m_i,(1,0,2)}(S) + \tfrac{1}{3} \beta_{m_j,(1,0,2)}(S).
\end{equation*}
Following the discussion in Section~\ref{sec:c1_smoothness}, the resulting spline is $\C{1}$-smooth at each vertex of $\triangle$ and also $\C{1}$-smooth across each edge of $\triangle$. Moreover, it is by construction $\C{1}$-smooth inside each triangle of $\triangle$ and consequently $\C{1}$-smooth everywhere on the domain $\Theta$.

We denote by $\mathbb{S}_3 \subset \splspace{3}{1}{\CTT}$ the vector space consisting of all splines that can be obtained by the described construction. It holds that $\dim(\mathbb{S}_3) = 3 n_v$ and $\poly_3 \subset \mathbb{S}_3$.

\subsection{Basis functions}
\label{sec:basis}

In the following we define basis functions for the \modified{spaces $\mathbb{S}_\ell$, $\ell = 1, 2, 3$. For a fixed $\ell$,} we associate with every vertex $\bfm{v}_i \in V$ three functions $B_{i,r}^{\ell}: \Theta \rightarrow \R$, $r = 0, 1, 2$\modified{, that form a basis of $\mathbb{S}_\ell$}.

The function $B_{i,r}^{\ell} \in \mathbb{S}_\ell$ is uniquely specified by choosing a linear polynomial $Q_{i,r}$ and setting $Q_i = Q_{i,r}$ and $Q_j = 0$, $j \neq i$, in the construction of splines described in Sections \ref{sec:s1}, \ref{sec:s2}, and \ref{sec:s3}. Independently of the choice of $Q_{i,r}$, this definition ensures that the support of $B_{i,r}^{\ell}$ is local. The construction in Section~\ref{sec:s1} implies that the support of $B_{i,r}^{1}$ is contained in the union of triangles of $\triangle$ that have a vertex at $\bfm{v}_i$. The supports of $B_{i,r} ^{2}$ and $B_{i,r}^{3}$ are, according to the constructions in Sections \ref{sec:s2} and \ref{sec:s3}, supplemented by the triangles of $\triangle$ that share an edge with a triangle that is in the support of $B_{i,r}^{1}$.

Any choice of linearly independent linear polynomials $Q_{i,r}$, $r = 0, 1, 2$, yields a set of linearly independent functions \modified{$B_{i,r}^\ell$}. Endorsing the ideas introduced in \cite{ps_dierckx_97}, we opt for polynomials determined by a barycentric coordinate system. Let $\convh{\bfm{q}_{i,0}, \bfm{q}_{i,1}, \bfm{q}_{i,2}} \subset \R^2$ be a triangle of a small area that contains $\bfm{v}_i$ and the points $\frac{5}{6} \bfm{v}_i + \frac{1}{6} \bfm{v}_j$ for every edge $\convh{\bfm{v}_i, \bfm{v}_j} \in E$. Such choice is motivated by the stability analysis in \cite{ps2_maes_04} and numerical tests made in Example~\ref{ex:vertex_triangles} (presented in Section~\ref{sec:examples}). Figure~\ref{fig:triangulation} (right) shows an example of such triangles that have minimal area.

\begin{figure}[!t]
\centering

\begin{subfigure}[b]{0.32\textwidth}
\centering
\includegraphics[width=\textwidth]{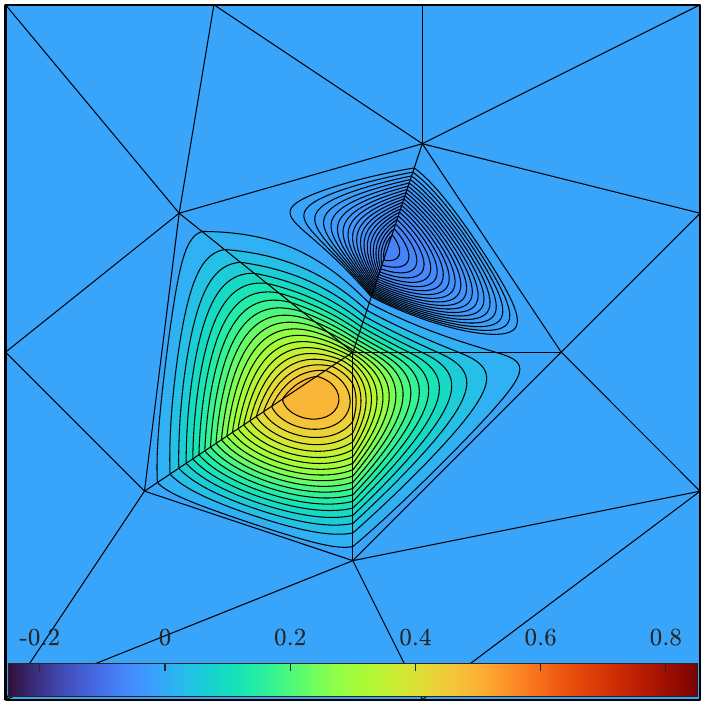}
\caption{Basis function $B_{i,0}^1$.}
\end{subfigure}
\hfill
\begin{subfigure}[b]{0.32\textwidth}
\centering
\includegraphics[width=\textwidth]{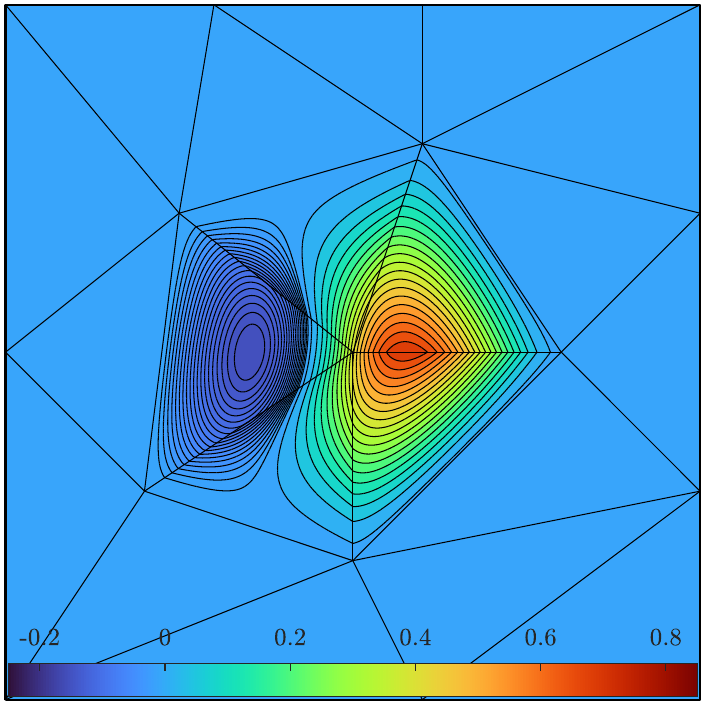}
\caption{Basis function $B_{i,1}^1$.}
\end{subfigure}
\hfill
\begin{subfigure}[b]{0.32\textwidth}
\centering
\includegraphics[width=\textwidth]{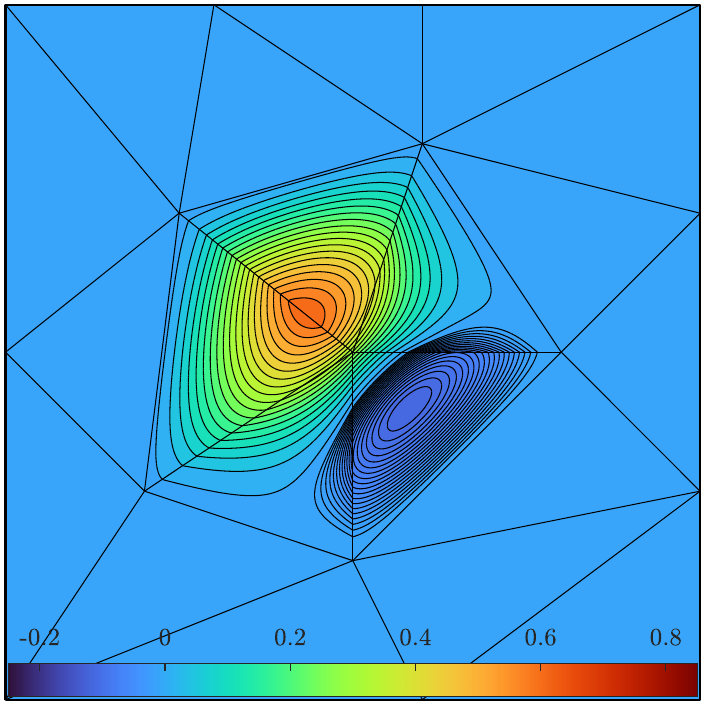}
\caption{Basis function $B_{i,2}^1$.}
\end{subfigure}

\begin{subfigure}[b]{0.32\textwidth}
\centering
\includegraphics[width=\textwidth]{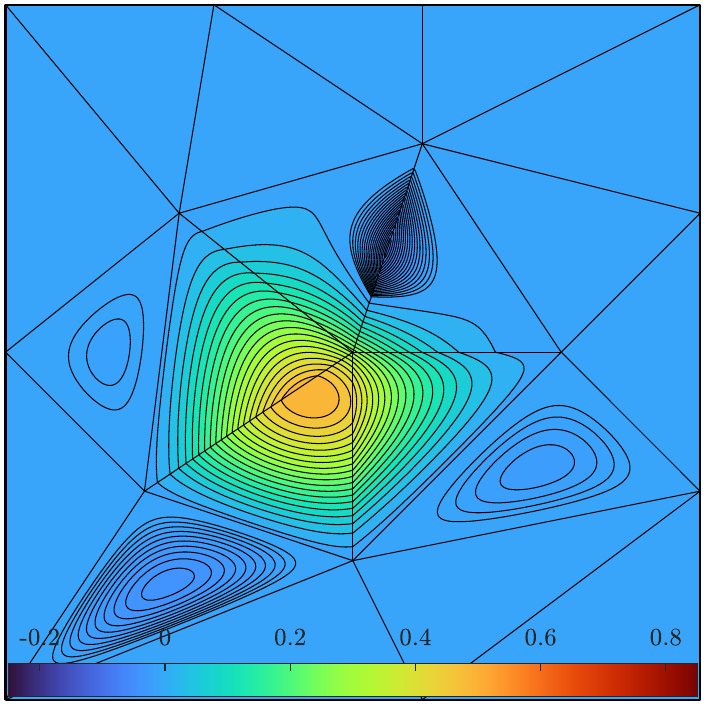}
\caption{Basis function $B_{i,0}^2$.}
\end{subfigure}
\hfill
\begin{subfigure}[b]{0.32\textwidth}
\centering
\includegraphics[width=\textwidth]{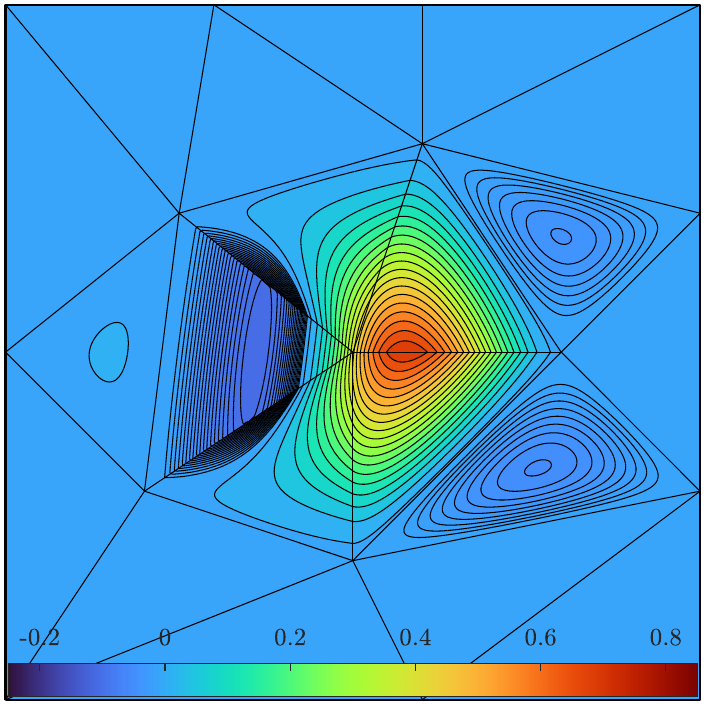}
\caption{Basis function $B_{i,1}^2$.}
\end{subfigure}
\hfill
\begin{subfigure}[b]{0.32\textwidth}
\centering
\includegraphics[width=\textwidth]{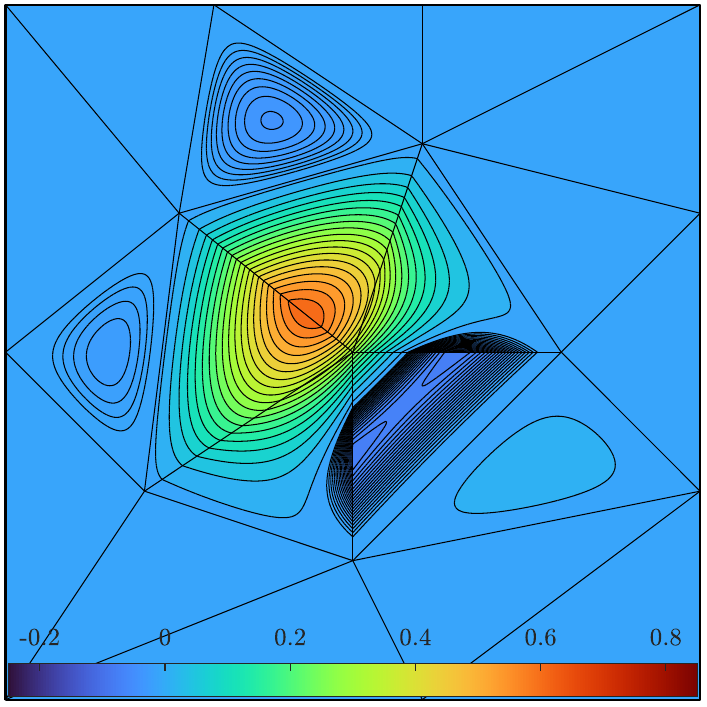}
\caption{Basis function $B_{i,2}^2$.}
\end{subfigure}

\begin{subfigure}[b]{0.32\textwidth}
\centering
\includegraphics[width=\textwidth]{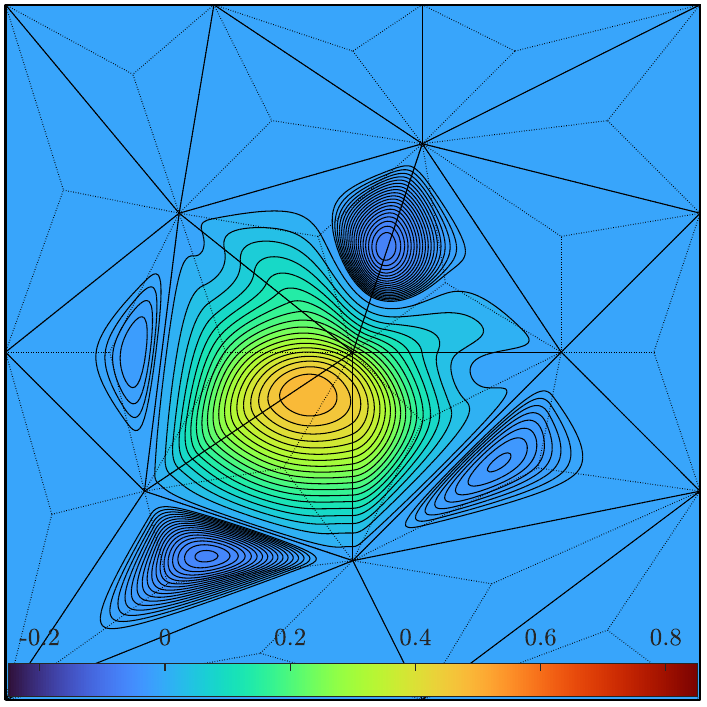}
\caption{Basis function $B_{i,0}^3$.}
\end{subfigure}
\hfill
\begin{subfigure}[b]{0.32\textwidth}
\centering
\includegraphics[width=\textwidth]{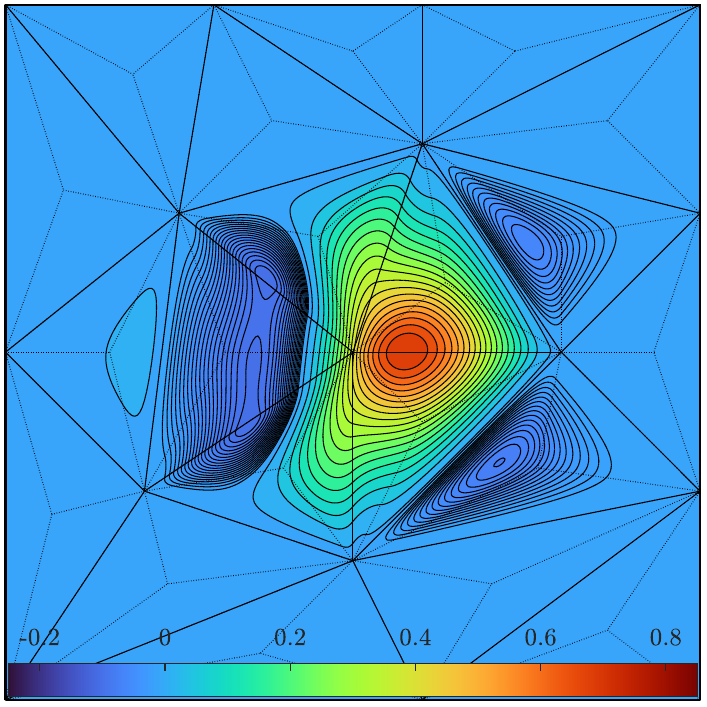}
\caption{Basis function $B_{i,1}^3$.}
\end{subfigure}
\hfill
\begin{subfigure}[b]{0.32\textwidth}
\centering
\includegraphics[width=\textwidth]{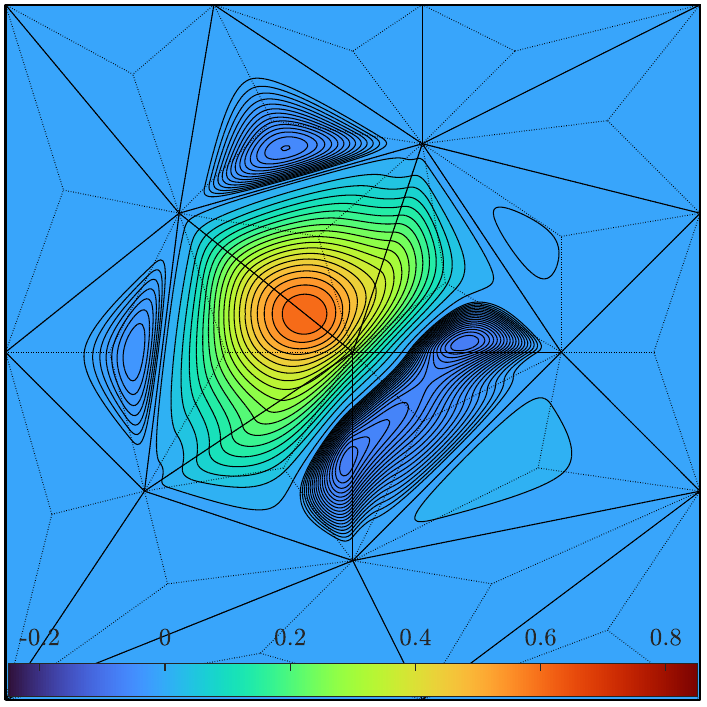}
\caption{Basis function $B_{i,2}^3$.}
\end{subfigure}

\caption{Contour plots of the basis functions $B_{i,r}^\ell$ for $r = 0, 1, 2$ and $\ell = 1, 2, 3$ that are associated with the vertex $\bfm{v}_i$ in the interior of the domain and defined based on the triangle $\convh{\bfm{q}_{i,0}, \bfm{q}_{i,1}, \bfm{q}_{i,2}}$ shown in Figure~\ref{fig:triangulation}.}
\label{fig:basis_interior_vertex}
\end{figure}

\begin{figure}[!t]
\centering

\begin{subfigure}[b]{0.32\textwidth}
\centering
\includegraphics[width=\textwidth]{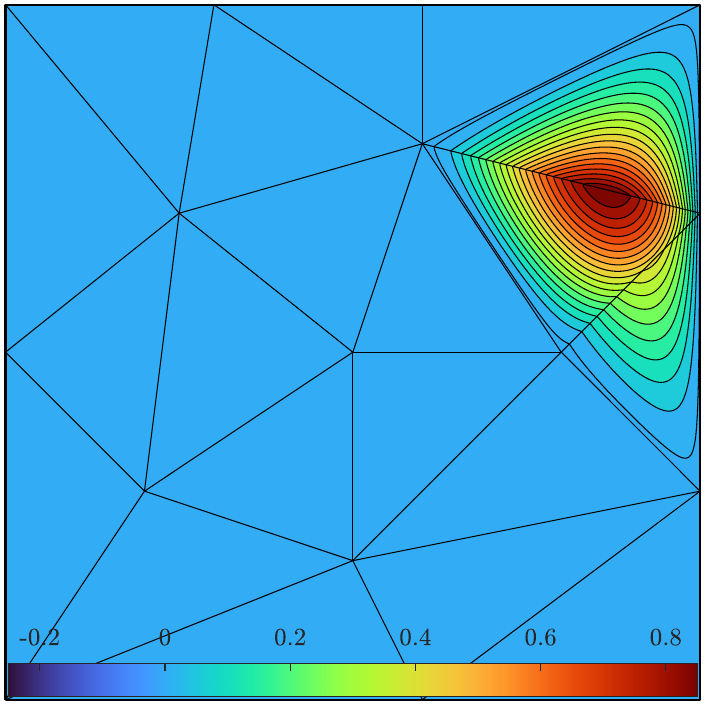}
\caption{Basis function $B_{j,0}^1$.}
\end{subfigure}
\hfill
\begin{subfigure}[b]{0.32\textwidth}
\centering
\includegraphics[width=\textwidth]{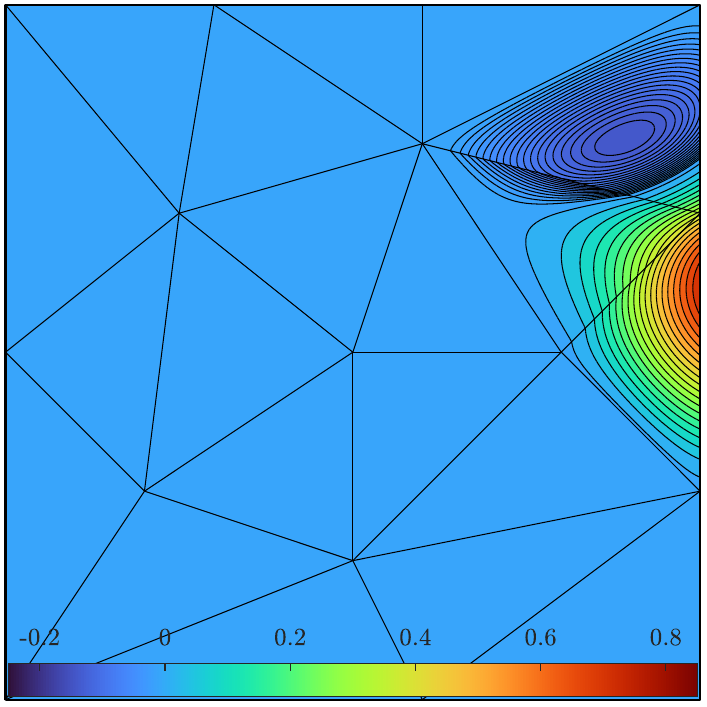}
\caption{Basis function $B_{j,1}^1$.}
\end{subfigure}
\hfill
\begin{subfigure}[b]{0.32\textwidth}
\centering
\includegraphics[width=\textwidth]{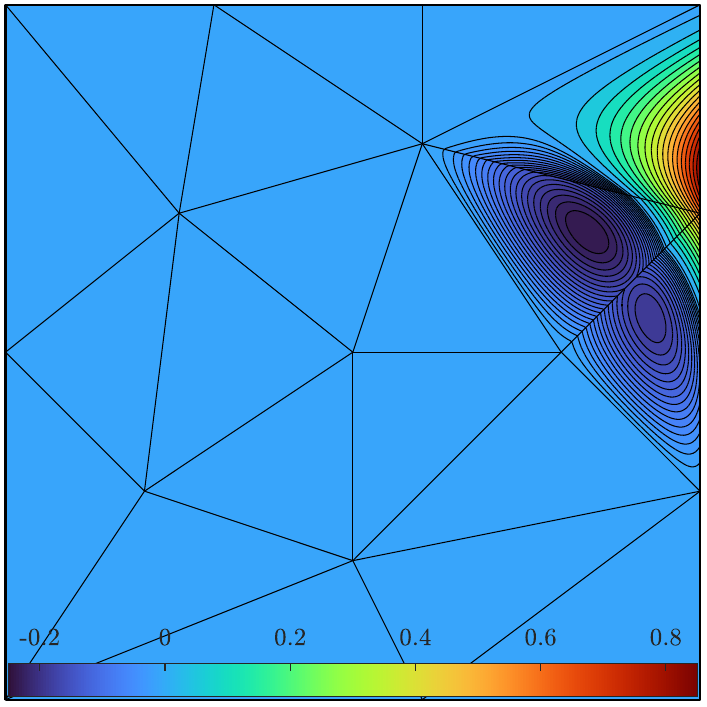}
\caption{Basis function $B_{j,2}^1$.}
\end{subfigure}

\begin{subfigure}[b]{0.32\textwidth}
\centering
\includegraphics[width=\textwidth]{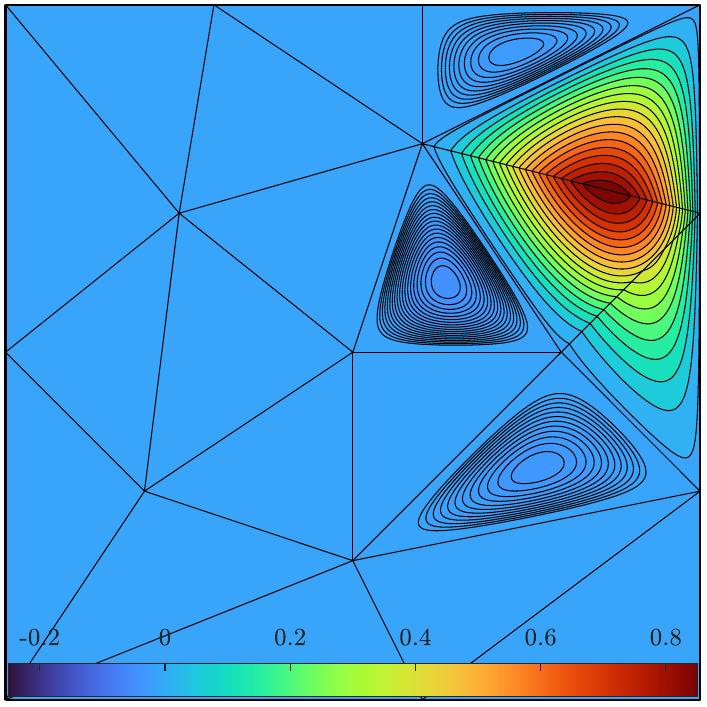}
\caption{Basis function $B_{j,0}^2$.}
\end{subfigure}
\hfill
\begin{subfigure}[b]{0.32\textwidth}
\centering
\includegraphics[width=\textwidth]{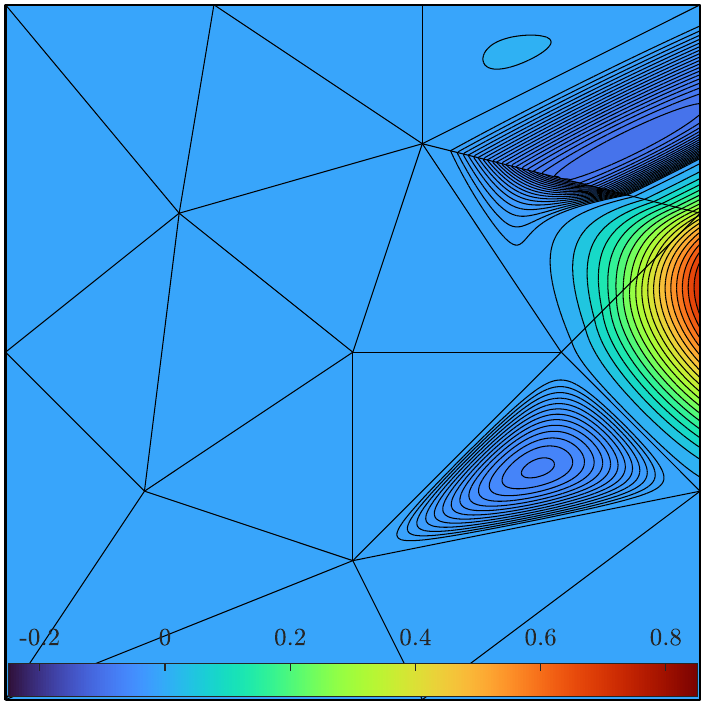}
\caption{Basis function $B_{j,1}^2$.}
\end{subfigure}
\hfill
\begin{subfigure}[b]{0.32\textwidth}
\centering
\includegraphics[width=\textwidth]{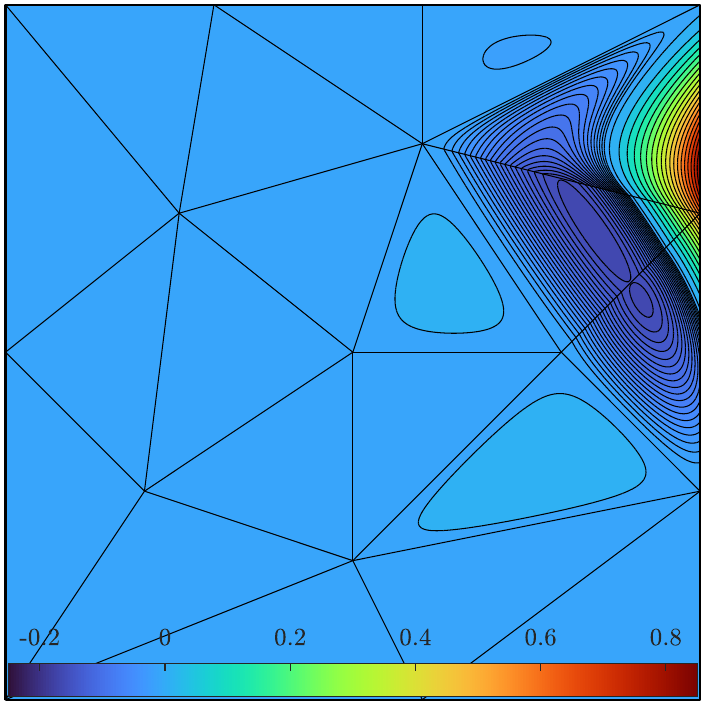}
\caption{Basis function $B_{j,2}^2$.}
\end{subfigure}

\begin{subfigure}[b]{0.32\textwidth}
\centering
\includegraphics[width=\textwidth]{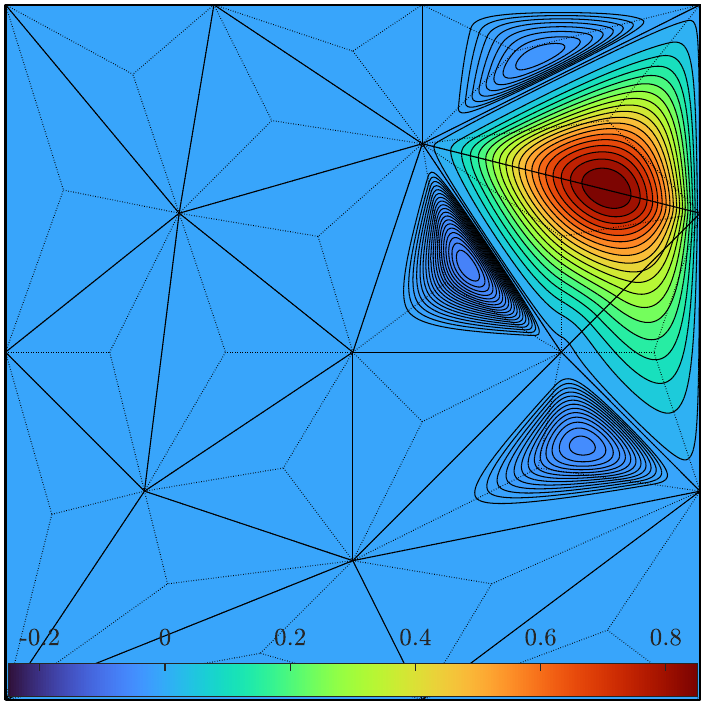}
\caption{Basis function $B_{j,0}^3$.}
\end{subfigure}
\hfill
\begin{subfigure}[b]{0.32\textwidth}
\centering
\includegraphics[width=\textwidth]{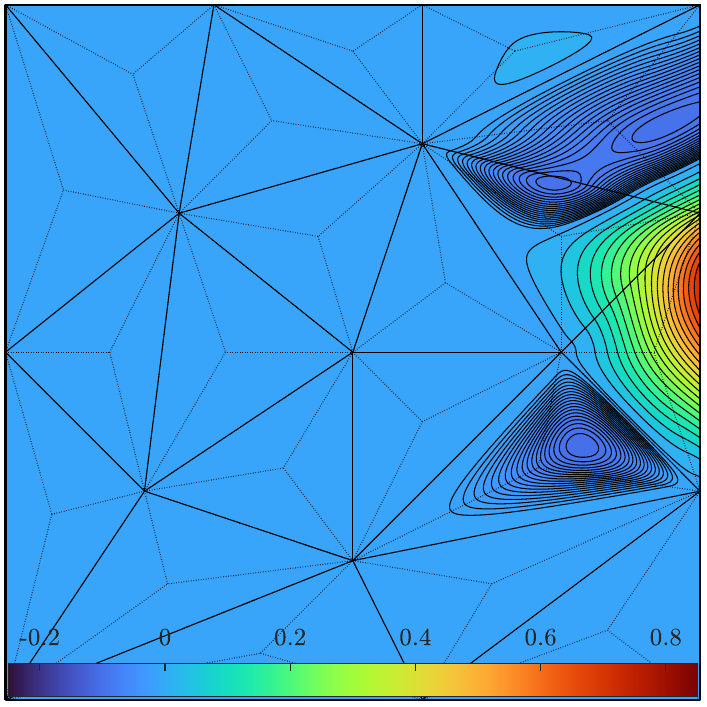}
\caption{Basis function $B_{j,1}^3$.}
\end{subfigure}
\hfill
\begin{subfigure}[b]{0.32\textwidth}
\centering
\includegraphics[width=\textwidth]{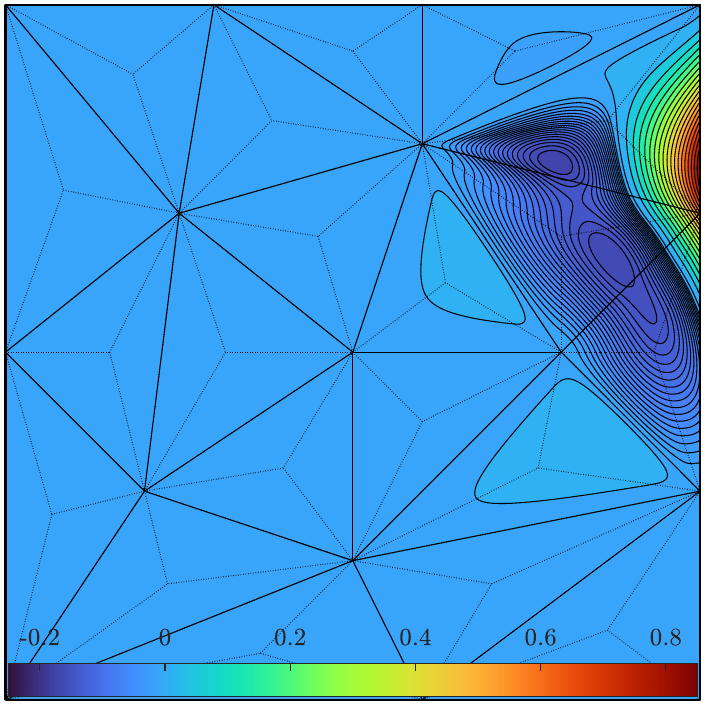}
\caption{Basis function $B_{j,2}^3$.}
\end{subfigure}

\caption{Contour plots of the basis functions $B_{j,r}^\ell$ for $r = 0, 1, 2$ and $\ell = 1, 2, 3$ that are associated with the vertex $\bfm{v}_j$ on the boundary of the domain and defined based on the triangle $\convh{\bfm{q}_{j,0}, \bfm{q}_{j,1}, \bfm{q}_{j,2}}$ shown in Figure~\ref{fig:triangulation}. Notice that the functions $B_{j,0}^\ell$ are zero on the boundary since the points $\bfm{q}_{j,1}$ and $\bfm{q}_{j,2}$ (marked in Figure~\ref{fig:triangulation}) lie on the boundary of the domain.}
\label{fig:basis_boundary_vertex}
\end{figure}

We define $(Q_{i,0}, Q_{i,1}, Q_{i,2})$ to be the barycentric coordinates with respect to $\convh{\bfm{q}_{i,0}, \bfm{q}_{i,1}, \bfm{q}_{i,2}}$. One obvious property of these linear polynomials is that $Q_{i,r}(\bfm{q}_{i,r}) = 1$ and $Q_{i,r}(\bfm{q}_{i,\rho}) = 0$ for every $\rho \in \set{0, 1, 2} \setminus \set{r}$. Examples of the resulting functions $B_{i,r}^\ell$ are shown in Figures \ref{fig:basis_interior_vertex} and \ref{fig:basis_boundary_vertex}.

\begin{remark}
In the construction of basis functions for the Powell--Sabin macro-element provided in \cite{ps_dierckx_97}, the triangle $\convh{\bfm{q}_{i,0}, \bfm{q}_{i,1}, \bfm{q}_{i,2}}$ associated with a vertex $\bfm{v}_i \in V$ is constrained to contain a certain set of points in order to ensure the nonnegativity of basis functions. This is not the case here as it is unclear how to achieve nonnegativity of the basis functions for $\mathbb{S}_2$ and $\mathbb{S}_3$. On the other hand, for the considered basis of $\mathbb{S}_1$, the nonnegativity can be easily enforced by requiring that $\convh{\bfm{q}_{i,0}, \bfm{q}_{i,1}, \bfm{q}_{i,2}}$ contains the point $\frac{2}{3} \bfm{v}_i + \frac{1}{3} \bfm{v}_j$ for every edge $\convh{\bfm{v}_i, \bfm{v}_j} \in E$ with a vertex at $\bfm{v}_i$ and the point $\frac{1}{2} \bfm{v}_i + \frac{1}{4} \bfm{v}_j + \frac{1}{4} \bfm{v}_k$ for every triangle $\convh{\bfm{v}_i, \bfm{v}_j, \bfm{v}_k} \in T$ with a vertex at $\bfm{v}_i$. This can be directly observed from \eqref{eq:s1_vertex_interpolation} and \eqref{eq:s1_central_value} by definition of the basis functions $B_{i,r}^1$, $r = 0, 1, 2$.
\end{remark}

\subsection{Spline representation}

A spline $S \in \mathbb{S}_\ell$, $\ell \in \set{1,2,3}$, can be uniquely represented as a linear combination of the basis functions $B_{i,r}^{\ell}$ for $\bfm{v}_i \in V$ and $r \in \set{0, 1, 2}$. To express the coefficients of this representation, we introduce the linear functionals
\begin{equation}
\label{eq:bfunc}
\beta_{i,r}^{\ell}(S) \modified{\coloneqq} S(\bfm{v}_i) + \grad{S}(\bfm{v}_i) \bfm{\cdot} (\bfm{q}_{i,r} - \bfm{v}_i).
\end{equation}
% where $\bfm{\cdot}$ denotes the inner product of two vectors.

Let $Q_i$ be the linear polynomial that by constructions in Sections \ref{sec:s1}, \ref{sec:s2}, \ref{sec:s3} specifies the value and the gradient of $S$ at $\bfm{v}_i$. From the discussion in Section~\ref{sec:c1_smoothness} we see that $\beta_{i,r}^\ell(S) = Q_i(\bfm{q}_{i,r})$.
\begin{comment}
Let $Q_i$ be the linear polynomial that by constructions in Sections \ref{sec:s1}, \ref{sec:s2}, \ref{sec:s3} specifies the value and the gradient of $S$ at $\bfm{v}_i$. Namely, $S(\bfm{v}_i) = Q_i(\bfm{v}_i)$ and
\begin{equation*}
\grad{S}(\bfm{v}_i) \bfm{\cdot} (\bfm{v}_j - \bfm{v}_i) = \grad{Q_i}(\bfm{v}_i) \bfm{\cdot} (\bfm{v}_j - \bfm{v}_i) = Q_i(\bfm{v}_j) - Q_i(\bfm{v}_i)
\end{equation*}
for every edge $\convh{\bfm{v}_i, \bfm{v}_j} \in E$. 
Suppose $\convh{\bfm{v}_i, \bfm{v}_j}$ and $\convh{\bfm{v}_i, \bfm{v}_k}$ are two distinct edges in $E$ with a vertex at $\bfm{v}_i$. There exist $\sigma_j, \sigma_k \in \R$ such that $\bfm{q}_{i,r} = (1 - \sigma_j - \sigma_k) \bfm{v}_i + \sigma_j \bfm{v}_j + \sigma_k \bfm{v}_k$. Then,
\begin{equation*}
\grad{S}(\bfm{v}_i) \bfm{\cdot} (\bfm{q}_{i,r} - \bfm{v}_i)
= \sigma_j \grad{S}(\bfm{v}_i) \bfm{\cdot} (\bfm{v}_j - \bfm{v}_i) + \sigma_k \grad{S}(\bfm{v}_i) \bfm{\cdot} (\bfm{v}_k - \bfm{v}_i)
= -(\sigma_j + \sigma_k) Q_i(\bfm{v}_i) + \sigma_j Q_i(\bfm{v}_j) + \sigma_k Q_i(\bfm{v}_k),
\end{equation*}
which implies that $\beta_{i,r}^\ell(S) = Q_i(\bfm{q}_{i,r})$.
\end{comment}
In particular, by definitions of the basis functions in Section~\ref{sec:basis}, we have
\begin{equation*}
\beta_{i,r}^{\ell}(B_{j,\rho}^{\ell}) =
\begin{cases}
1 & \text{if $i = j$ and $r = \rho$,} \\
0 & \text{otherwise,}
\end{cases}
\end{equation*}
and hence
\begin{equation}
\label{eq:brep}
S = \sum_{i = 1}^{n_v} \sum_{r = 0}^2 \beta_{i,r}^{\ell}(S) B_{i,r}^{\ell}.
\end{equation}
This result has some interesting consequences.

As the constant $(x, y) \mapsto 1$ is an element of $\mathbb{S}_\ell$, \eqref{eq:bfunc} and \eqref{eq:brep} imply
\begin{equation*}
1 = \sum_{i = 1}^{n_v} \sum_{r = 0}^2 B_{i,r}^{\ell}(x,y)
\end{equation*}
for every $(x,y) \in \Theta$, i.e., the basis functions form a partition of unity. Moreover, considering the linear polynomials $(x,y) \mapsto x$ and $(x,y) \mapsto y$, it follows by the same line of arguments that
\begin{equation*}
(x,y) = \sum_{i = 1}^{n_v} \sum_{r = 0}^2 \bfm{q}_{i,r} B_{i,r}^{\ell}(x,y)
\end{equation*}
for every $(x,y) \in \Theta$, \ie{}, the points $\bfm{q}_{i,r}$, $r = 0, 1, 2$, that define the triangle associated with $\bfm{v}_i \in V$ are the Greville points of the basis functions $B_{i,r}^\ell$.

\section{Numerical comparison}
\label{sec:examples}

This section provides a number of numerical examples with the aim to demonstrate the behavior of the cubic splines presented in Section~\ref{sec:cubic_splines}. We consider the problems of finding the best (discrete) $L^2$ approximation to a function and solving the Poisson boundary value problem via the isoparametric finite element method and the immersed penalized boundary method.

\subsection{Preliminaries}

In what follows we deal with methods that require the representation of a spline $S \in \mathbb{S}_\ell$, $\ell \in \set{1,2,3}$, in terms of basis functions. For this purpose we use the bases introduced in Section~\ref{sec:basis}. To make the notation simpler, we denote the functions $B_{i,r}^\ell$, $i = 1, 2, \ldots, n_v$, $r = 0, 1, 2$, for a fixed $\ell$ by $S_k$ with the indices $k$ ranging from $1$ to $N = 3 n_v$. The spline $S$ is expressed in the form
\begin{equation}
\label{eq:srep}
S = \sum_{k = 1}^N c_k S_k
\end{equation}
for coefficients $c_k \in \R$. Additionally, we assume that the basis functions are enumerated in such a way that $S_1, S_2, \ldots, S_{n}$ for $n < N$ are zero on the boundary of $\Theta$, and $S_{n+1}, S_{n+2}, \ldots, S_N$ are nonzero on the boundary of $\Theta$. This arrangement is important for the finite element method described in Section~\ref{sec:boundary_value_problems}.

\begin{example}
\label{ex:triangulation}
Let $\Theta$ be the unit square partitioned by the triangulation $\triangle$ shown in Figure~\ref{fig:triangulation}. The triangulation has $16$ vertices, hence $N = 48$. The basis functions are determined by the shaded triangles, one for each vertex. For each vertex on the boundary of $\Theta$ that is shared by two collinear edges, \ie{}, a boundary vertex that is not a corner of the square, the triangle is positioned in a way that two of its vertices lie on the boundary. This means that of the three basis functions associated with the vertex only two are nonzero on the boundary (see Figure~\ref{fig:basis_boundary_vertex}), which implies that $n = 24$. An important aspect of such choice is that it in general ensures that the restrictions of $S_{n+1}, S_{n+2}, \ldots, S_N$ to $\partial \Theta$ are linearly independent (see \cite{ct_groselj_22} for a more detailed discussion).
\end{example}

\begin{figure}[!t]
\centering

\begin{subfigure}[b]{0.32\textwidth}
\centering
\includegraphics[width=\textwidth]{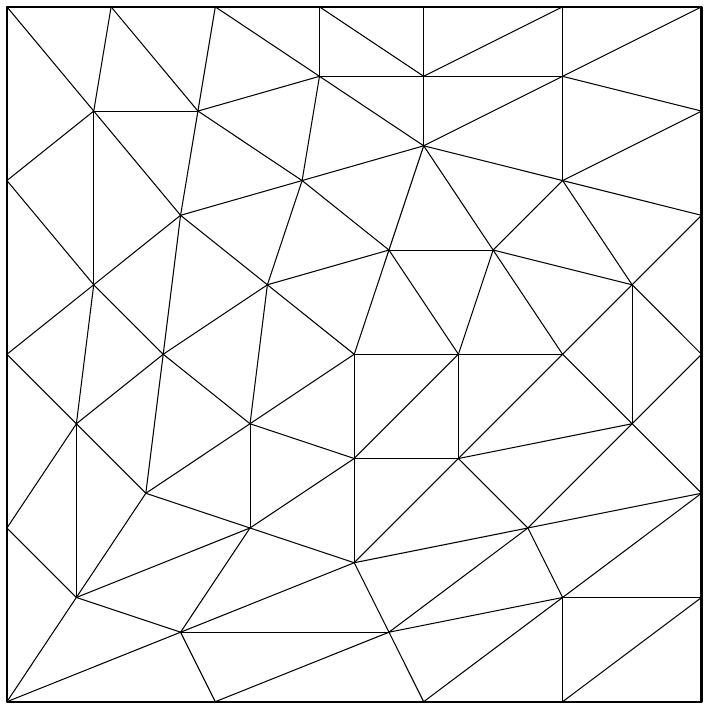}
\end{subfigure}
\hfill
\begin{subfigure}[b]{0.32\textwidth}
\centering
\includegraphics[width=\textwidth]{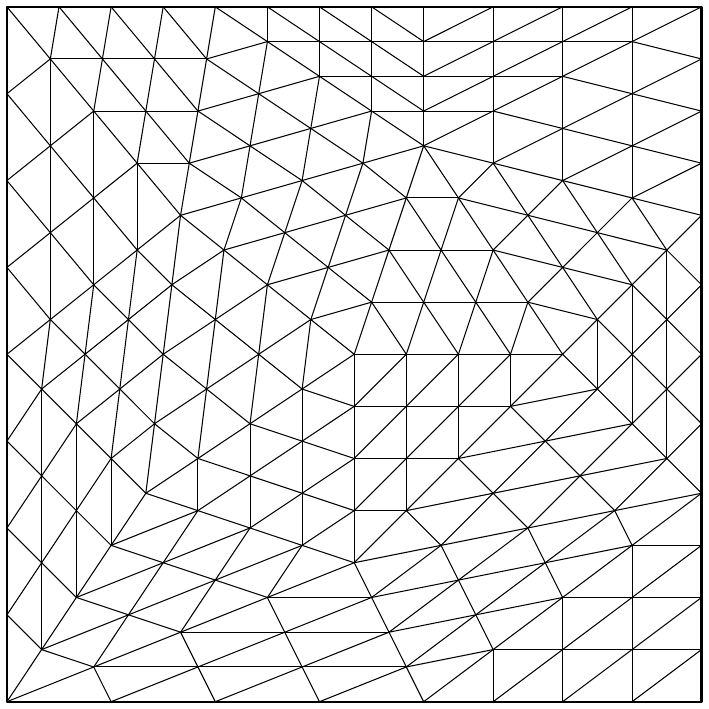}
\end{subfigure}
\hfill
\begin{subfigure}[b]{0.32\textwidth}
\centering
\includegraphics[width=\textwidth]{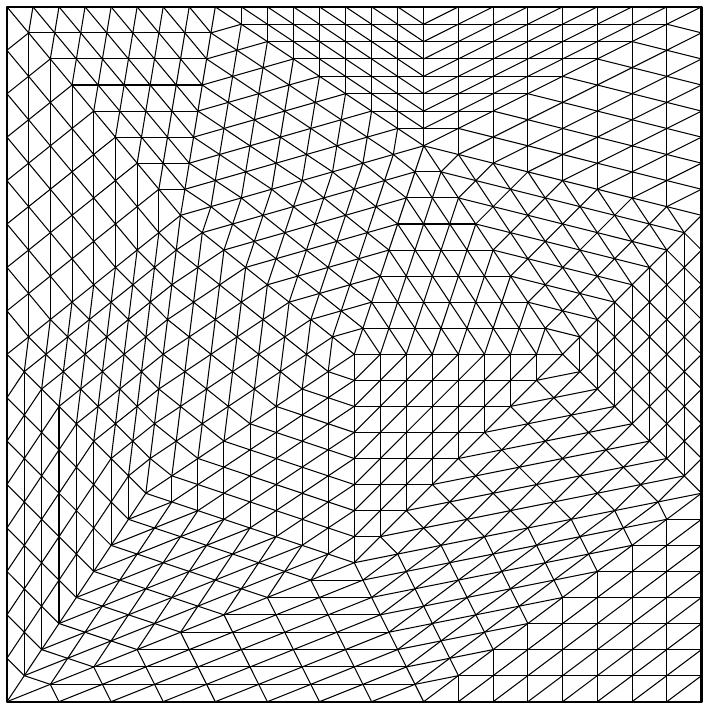}
\end{subfigure}

\caption{A sequence of triangulations obtained by subsequently refining the triangulation of the unit square shown in Figure~\ref{fig:triangulation}.}
\label{fig:triangulation_sequence}
\end{figure}

\begin{example}
\label{ex:triangulation_sequence}
Figure~\ref{fig:triangulation_sequence} shows three subsequent refinements of the triangulation presented in Example~\ref{ex:triangulation} and together with it forms a sequence of four triangulations. In each step we refine the triangulation by subdividing each triangle into four smaller triangles determined by the vertices and midpoints of the edges. We also adjust and supplement the set of triangles associated with the vertices so that the basis functions on the refined triangulations satisfy the properties discussed in Example~\ref{ex:triangulation}.
\end{example}

In the process of refinement described in Example~\ref{ex:triangulation_sequence}, the length $h$ of the longest edge in the triangulation is halved, which allows us to numerically inspect the convergence rates of the methods. With respect to the polynomial precision properties, it is expected that the methods based on $\mathbb{S}_1$ converge with rate at most $3$, and the methods based on $\mathbb{S}_2$ and $\mathbb{S}_3$ converge with rate at most $4$, \ie{}, the approximation error in the best case behaves as $\mc{O}(h^3)$ and $\mc{O}(h^4)$, respectively. In the numerical tests conducted over the triangulations from Examples \ref{ex:triangulation} and \ref{ex:triangulation_sequence}, the approximation error is measured as the maximal absolute difference between exact solution and approximation computed on the $401 \times 401$ uniform grid of points in the unit square.

The methods described in Sections \ref{sec:function_approximations} and \ref{sec:boundary_value_problems} reduce to solving a system of linear equations $\bfm{A} \, \bfm{c} = \bfm{b}$. The matrix $\bfm{A}$ is of size either $N \times N$ or $n \times n$, and the vectors $\bfm{c}$ and $\bfm{b}$ are of length either $N$ or $n$. In this respect the elements of $\bfm{c}$ are either all or the first $n$ coefficients $c_k$ in \eqref{eq:srep} collected increasingly by $k$. The elements of $\bfm{A}$ are denoted by either $a_{j,k}^N$ or $a_{j,k}^n$, and the elements of $\bfm{b}$ are denoted by either $b_j^N$ or $b_j^n$ with the indices $j$ and $k$ ranging from $1$ to $N$ or $n$, respectively. We measure stability of the spline representation by the $2$-norm condition number (cond) of the matrix $\bfm{A}$.

\begin{figure}
\centering

\begin{subfigure}[b]{0.32\textwidth}
\centering
\includegraphics[width=\textwidth]{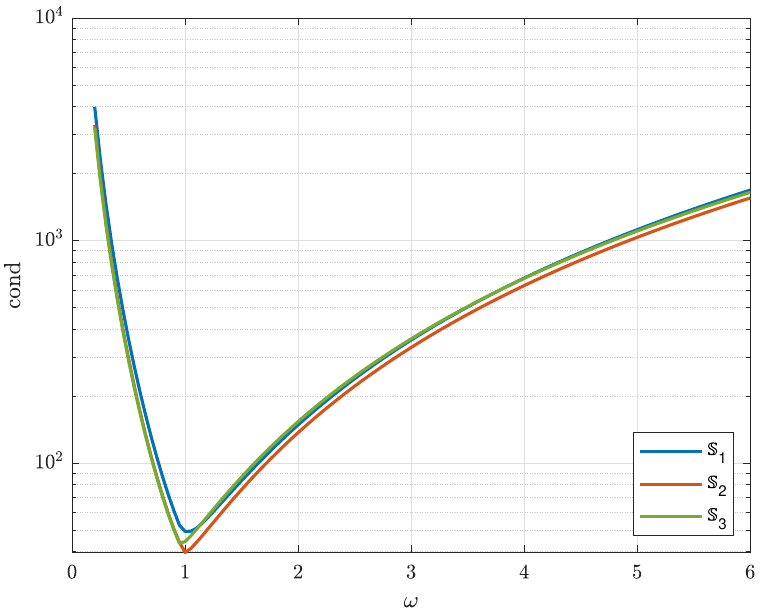}
\caption{Best $L^2$ approximation.}
\end{subfigure}
\hfill
\begin{subfigure}[b]{0.32\textwidth}
\centering
\includegraphics[width=\textwidth]{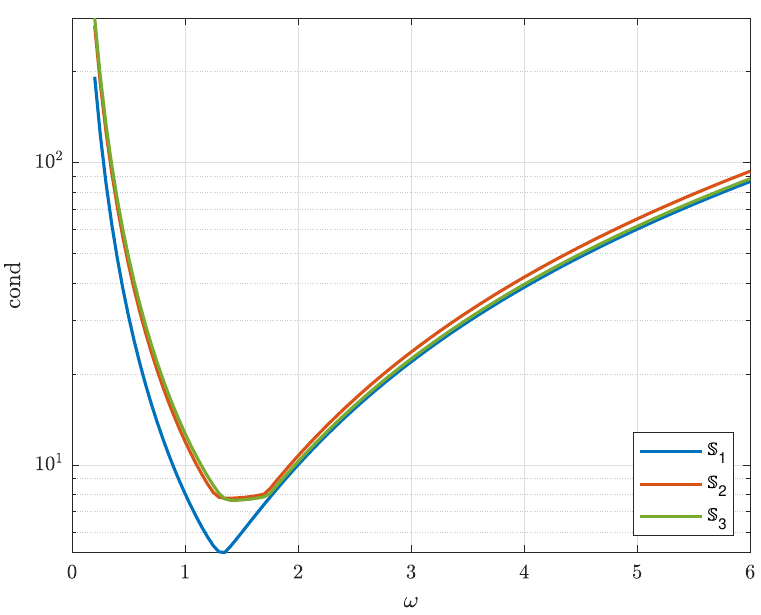}
\caption{Finite element method.}
\end{subfigure}
\hfill
\begin{subfigure}[b]{0.32\textwidth}
\centering
\includegraphics[width=\textwidth]{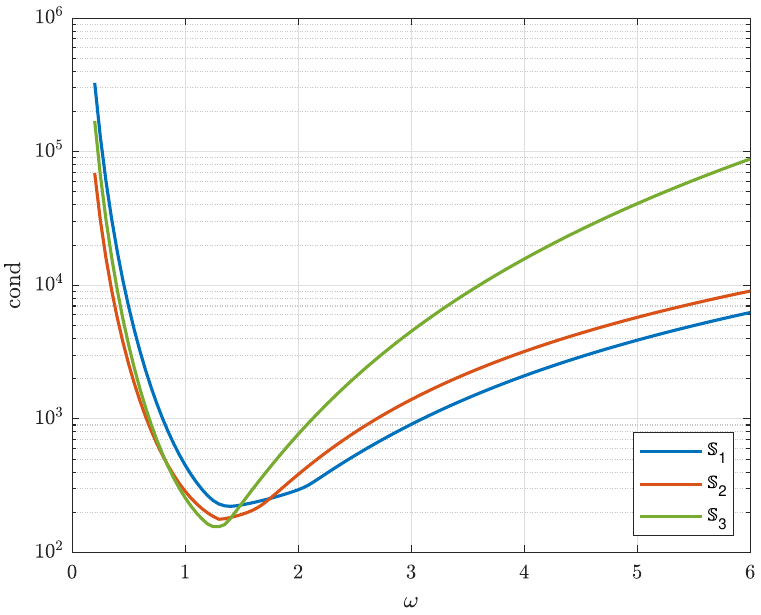}
\caption{Immersed penalized boundary method.}
\end{subfigure}

\caption{The condition numbers of matrices arising in different approximation methods with respect to the scaling factor $\omega$ that determines the size of triangles used in the construction of the basis functions for $\mathbb{S}_\ell$, $\ell = 1, 2, 3$, as described in Example~\ref{ex:vertex_triangles}.}
\label{fig:vertex_triangles}
\end{figure}

\begin{example}
\label{ex:vertex_triangles}
For the triangulation discussed in Example~\ref{ex:triangulation} we inspect the condition numbers of the matrices $\bfm{A}$ arising in the best $L^2$ approximation (Example~\ref{ex:l2_franke}), finite element method (Example~\ref{ex:fem}), and immersed penalized boundary method (Example~\ref{ex:ipbm}) with respect to the size of the triangles that define the basis functions. Suppose $\convh{\bfm{q}_{i,0}, \bfm{q}_{i,1}, \bfm{q}_{i,2}}$ are the triangles associated with the vertices $\bfm{v}_i \in V$ depicted in Figure~\ref{fig:triangulation} (right), and let
\begin{equation*}
\bfm{q}_{i,r}^\omega = \bfm{v}_i + \omega (\bfm{q}_{i,r} - \bfm{v}_i), \quad r = 0, 1, 2,
\end{equation*}
for a factor $\omega > 0$. Figure~\ref{fig:vertex_triangles} depicts the condition numbers in dependence of $\omega \in [0.2, 6]$ for the considered basis of $\mathbb{S}_\ell$, $\ell \in \set{1, 2, 3}$, determined by the triangles $\convh{\bfm{q}_{i,0}^\omega, \bfm{q}_{i,1}^\omega, \bfm{q}_{i,2}^\omega}$. The results indicate that $\omega = 1$ ensures near-optimal stability.
\end{example}

The expressions for the elements of $\bfm{A}$ and $\bfm{b}$ typically include integrals over a domain. If $S$ is a spline over a triangulation $\triangle$ of the domain $\Theta$ that is defined piecewise on each triangle $\bfm{t}_m$, $m \in \set{1, 2, \ldots, n_t}$, of $\triangle$ by a polynomial $P_m$, the integral of $S$ over the domain is always considered as the sum of integrals of $P_m$ over the triangles of $\triangle$. More precisely,
\begin{equation*}
\int_\Theta S \, \mathrm{d}\Theta = \sum_{m = 1}^{n_t} \bracket{\int_{\bfm{t}_m} P_m \, \mathrm{d}\bfm{t}_m}.
\end{equation*}
This way the integral exists even if $S$ is not well defined on the edges of $\triangle$.

In the following we briefly discuss some common methods for function approximation and solving boundary value problems and provide several numerical examples. We comment the results in Section~\ref{sec:conclusion}.

\subsection{Function approximations}
\label{sec:function_approximations}

Let $\Theta$ be a polygonal domain and $\triangle$ its triangulation. The best $L^2$ approximation $S \in \mathbb{S}_\ell$, $\ell \in \set{1,2,3}$, of a given function $F \in L^2(\Theta)$, \ie{}, the solution to the minimization problem
\begin{equation*}
\underset{S \in \mathbb{S}_\ell}{\mathrm{arg\,min}} \int_\Theta \bracket{F-S}^2 \, \mathrm{d}\Theta,
\end{equation*}
is determined by the conditions
\begin{equation}
\label{eq:l2_approximation_conditions}
\int_\Theta \bracket{F - S} \, S_j \, \mathrm{d}\Theta = 0, \quad j = 1, 2, \ldots, N.
\end{equation}
If we represent $S$ in the form \eqref{eq:srep}, the computation of the unknown coefficients $c_k$ amounts to solving the system $\bfm{A} \, \bfm{c} = \bfm{b}$ determined by
\begin{equation*}
a_{j,k}^N = \int_\Theta S_k \, S_j \, \mathrm{d}\Theta, \quad b_j^N = \int_\Theta F \, S_j \, \mathrm{d}\Theta.
\end{equation*}
The matrix $\bfm{A}$ is the Gram matrix associated with the functions $S_1, S_2, \ldots, S_N$.

\begin{figure}[!t]
\centering

\begin{subfigure}[b]{0.49\textwidth}
\centering
\includegraphics[width=\textwidth]{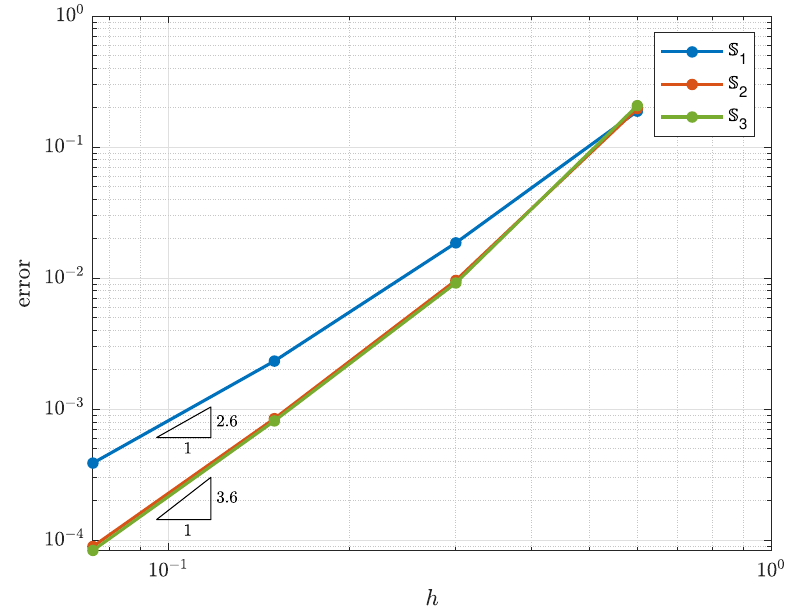}
\caption{Errors for Example~\ref{ex:l2_franke}.}
\end{subfigure}
\hfill
\begin{subfigure}[b]{0.49\textwidth}
\centering
\includegraphics[width=\textwidth]{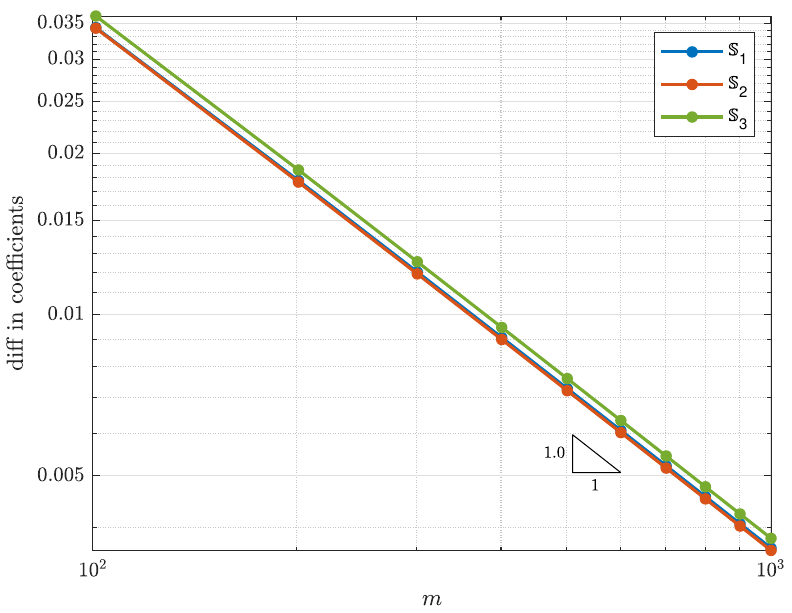}
\caption{Convergence discussed in Example~\ref{ex:dl2_franke}.}
\end{subfigure}

\caption{Errors of the approximations from $\mathbb{S}_\ell$, $\ell = 1, 2, 3$, for the problem discussed in Example~\ref{ex:l2_franke} (left) and the convergence of the coefficient vectors corresponding to the discrete best $L^2$ approximation when the number of points is increased as described in Example~\ref{ex:dl2_franke} (right).}
\label{fig:l2_franke_errors}
\end{figure}

\begin{figure}[!t]
\centering

\begin{subfigure}[b]{0.24\textwidth}
\centering
\includegraphics[width=\textwidth]{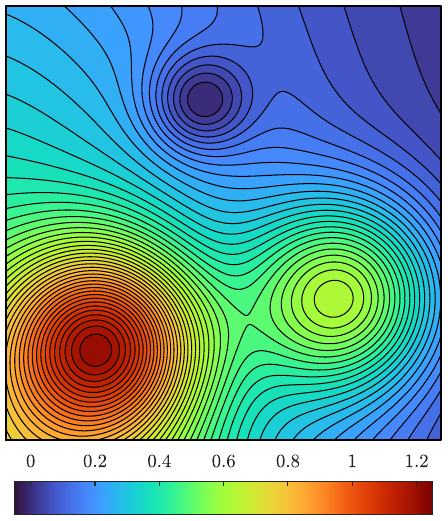}
\caption{Franke's test function.}
\end{subfigure}
\hfill
\begin{subfigure}[b]{0.24\textwidth}
\centering
\includegraphics[width=\textwidth]{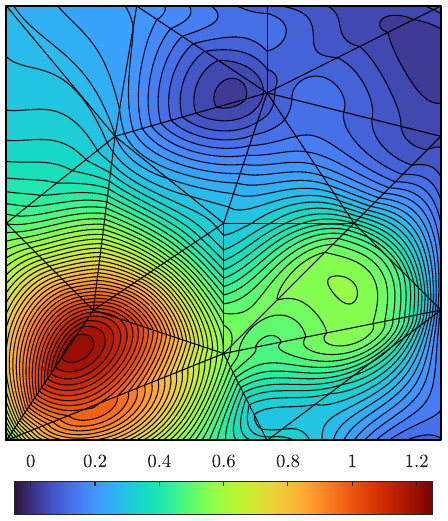}
\caption{Best $L^2$ approximation in $\mathbb{S}_1$.}
\end{subfigure}
\hfill
\begin{subfigure}[b]{0.24\textwidth}
\centering
\includegraphics[width=\textwidth]{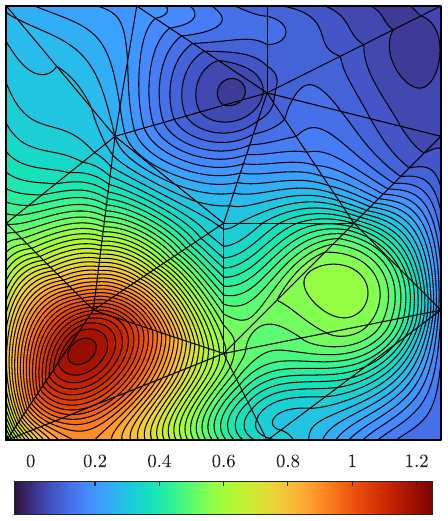}
\caption{Best $L^2$ approximation in $\mathbb{S}_2$.}
\end{subfigure}
\hfill
\begin{subfigure}[b]{0.24\textwidth}
\centering
\includegraphics[width=\textwidth]{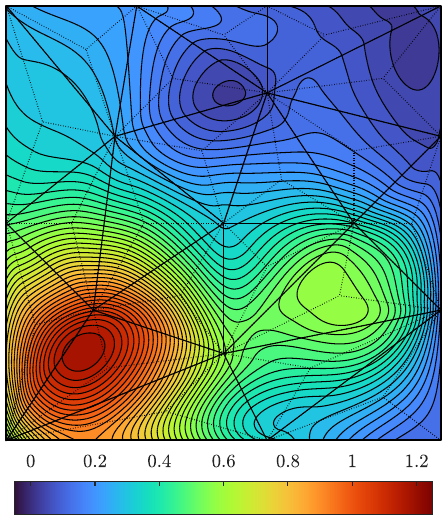}
\caption{Best $L^2$ approximation in $\mathbb{S}_3$.}
\end{subfigure}

\caption{Contour plots of the Franke's test function and its best $L^2$ approximations in $\mathbb{S}_\ell$, $\ell = 1, 2, 3$.}
\label{fig:l2_franke}
\end{figure}

\begin{example}
\label{ex:l2_franke}
We compute the best $L^2$ approximation in $\mathbb{S}_\ell$, $\ell = 1, 2, 3$, to the Franke's test function
\begin{equation}
\label{eq:franke_fun}
\begin{aligned}
F(x,y) &= 0.75 \, \mathrm{exp}\bracket{-\tfrac{1}{4}(9x - 2)^2 - \tfrac{1}{4}(9y - 2)^2}
+ 0.75 \, \mathrm{exp}\bracket{-\tfrac{1}{49} (9x + 1)^2 - \tfrac{1}{10} (9y + 1)} \\[0.2em]
& \quad + 0.5 \, \mathrm{exp}\bracket{-\tfrac{1}{4} (9x - 7)^2 - \tfrac{1}{4} (9y - 3)^2}
- 0.2 \, \mathrm{exp}\bracket{-(9x - 4)^2 - (9y - 7)^2}
\end{aligned}
\end{equation}
on the unit square partitioned by the sequence of triangulations discussed in Example~\ref{ex:triangulation_sequence}. The results are provided in Table~\ref{tab:l2_franke}. Figure~\ref{fig:l2_franke_errors} (left) depicts the errors with respect to the length of the longest edge of the underlying triangulation. Figure~\ref{fig:l2_franke} shows the contour plots of the approximations over the initial triangulation.
\end{example}

\begin{table}[!t]
\centering
\begin{tabular}{|c|cc|cc|cc|}
\hline
& \multicolumn{2}{c|}{$\mathbb{S}_1$} & \multicolumn{2}{c|}{$\mathbb{S}_2$} & \multicolumn{2}{c|}{$\mathbb{S}_3$} \\
$N$ & error & cond & error & cond & error & cond \\ \hline 
48 & 1.88e-01 & 4.94e+01  & 1.97e-01 & 3.97e+01  & 2.07e-01 & 4.45e+01 \\
153 & 1.86e-02 & 4.58e+01  & 9.62e-03 & 4.12e+01  & 9.20e-03 & 4.25e+01 \\
543 & 2.33e-03 & 4.69e+01  & 8.49e-04 & 4.66e+01  & 8.15e-04 & 5.08e+01 \\
2043 & 3.89e-04 & 4.95e+01  & 8.99e-05 & 5.04e+01  & 8.39e-05 & 5.71e+01 \\
\hline
\end{tabular}
\caption{Approximation errors and condition numbers for the problem described in Example~\ref{ex:l2_franke}.}
\label{tab:l2_franke}
\end{table}

In the discrete version of the best $L^2$ approximation we use given values $F(\bfm{p}_d)$, $d = 1, 2, \ldots, M$, at $M \geq N$ points $\bfm{p}_d \in \Theta$ and search for
\begin{equation*}
\underset{S \in \mathbb{S}_\ell}{\mathrm{arg\,min}} \sum_{d = 1}^M \bracket{F(\bfm{p}_d) - S(\bfm{p}_d)}^2.
\end{equation*}
To this end, we replace the conditions in \eqref{eq:l2_approximation_conditions} by
\begin{equation*}
\sum_{d = 1}^M \bracket{F(\bfm{p}_d) - S(\bfm{p}_d)} S_j(\bfm{p}_d) = 0, \quad j = 1, 2, \ldots, N.
\end{equation*}
This yields the system $\bfm{A} \, \bfm{c} = \bfm{b}$ determined by
\begin{equation*}
a_{j,k}^N = \sum_{d = 1}^M S_k(\bfm{p}_d) \, S_j(\bfm{p}_d), \quad b_j^N = \sum_{d = 1}^M F(\bfm{p}_d) \, S_j(\bfm{p}_d).
\end{equation*}
Alternatively, the vector of coefficients $\bfm{c}$ can be expressed as the best least squares solution to the overdetermined system $\bfm{S} \, \bfm{c} = \bfm{f}$, where $\bfm{S} \in \R^{M \times N}$ is the collocation matrix determined by element $S_k(\bfm{p}_d)$ at position $(d,k)$, and $\bfm{f} \in \R^M$ is the vector with element $F(\bfm{p}_d)$ at position $d$. 

\begin{example}
\label{ex:dl2_franke}
We repeat Example~\ref{ex:l2_franke} for the discrete best $L^2$ approximation based on the $201 \times 201$ uniform grid of points in the unit square. Table~\ref{tab:dl2_franke} reveals that results are similar to the ones observed in Example~\ref{ex:l2_franke}. For the triangulation discussed in Example~\ref{ex:triangulation}, Figure~\ref{fig:l2_franke_errors} (right) shows linear decay in $2$-norm of the difference between the coefficient vector corresponding to the best $L^2$ approximation and the coefficient vector corresponding to the discrete best $L^2$ approximation on the $m \times m$ uniform grid of points ($M = m^2$) for $m \in \set{101, 201, \ldots, 1001}$. 
\end{example}

\begin{table}[!t]
\centering
\begin{tabular}{|c|cc|cc|cc|}
\hline
& \multicolumn{2}{c|}{$\mathbb{S}_1$} & \multicolumn{2}{c|}{$\mathbb{S}_2$} & \multicolumn{2}{c|}{$\mathbb{S}_3$} \\
$N$ & error & cond & error & cond & error & cond \\ \hline 
48 & 1.76e-01 & 4.63e+01  & 1.85e-01 & 3.69e+01  & 1.95e-01 & 4.15e+01 \\
153 & 1.87e-02 & 3.91e+01  & 9.59e-03 & 3.41e+01  & 9.20e-03 & 3.63e+01 \\
543 & 2.33e-03 & 3.87e+01  & 8.48e-04 & 3.45e+01  & 8.18e-04 & 3.81e+01 \\
2043 & 3.90e-04 & 3.90e+01  & 9.00e-05 & 3.55e+01  & 8.38e-05 & 4.21e+01 \\
\hline
\end{tabular}
\caption{Approximation errors and condition numbers for the problem described in Example~\ref{ex:dl2_franke}.}
\label{tab:dl2_franke}
\end{table}

Since scattered function data to be approximated is often subjected to noise, it is common to use a penalty term in the minimization of the discrete $L^2$ norm to prevent high oscillation of the approximant. Following \cite{schumaker_15}, we can search for
\begin{equation*}
\underset{S \in \mathbb{S}_\ell}{\mathrm{arg\,min}} \bracket{
\sum_{d = 1}^M \bracket{F(\bfm{p}_d) - S(\bfm{p}_d)}^2 +
\lambda \int_{\Theta} \bfm{e} \bfm{\cdot} (\bfm{H}_S \bfm{\odot} \bfm{H}_S) \bfm{\cdot} \bfm{e} \ \mathrm{d}\Theta
},
\end{equation*}
where $\lambda \geq 0$ is the penalty parameter, $\bfm{e} = (1, 1)$, $\bfm{H}_S$ is the Hessian matrix of $S$, and $\bfm{\odot}$ denotes the Hadamard product. The solution to the penalized discrete best $L^2$ approximation problem can be computed by solving the system $\bfm{A} \, \bfm{c} = \bfm{b}$ for
\begin{equation*}
\bfm{A} = \bfm{S}^T \bfm{S} + \lambda \, \bfm{H}, \quad
\bfm{b} = \bfm{S}^T \bfm{f},
\end{equation*}
where $\bfm{H} \in \R^{N \times N}$ is determined by the element $\int_\Theta \bfm{e} \bfm{\cdot} (\bfm{H}_{S_k} \bfm{\odot} \bfm{H}_{S_j}) \bfm{\cdot} \bfm{e} \, \mathrm{d}\Theta$ at position $(j,k)$.

\begin{figure}[!t]
\centering

\begin{subfigure}[b]{0.49\textwidth}
\centering
\includegraphics[width=\textwidth]{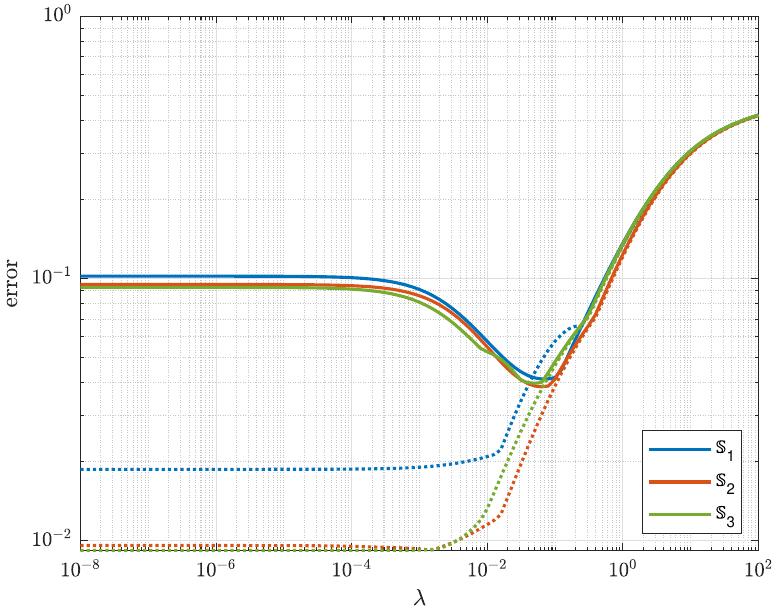}
\end{subfigure}
\hfill
\begin{subfigure}[b]{0.49\textwidth}
\centering
\includegraphics[width=\textwidth]{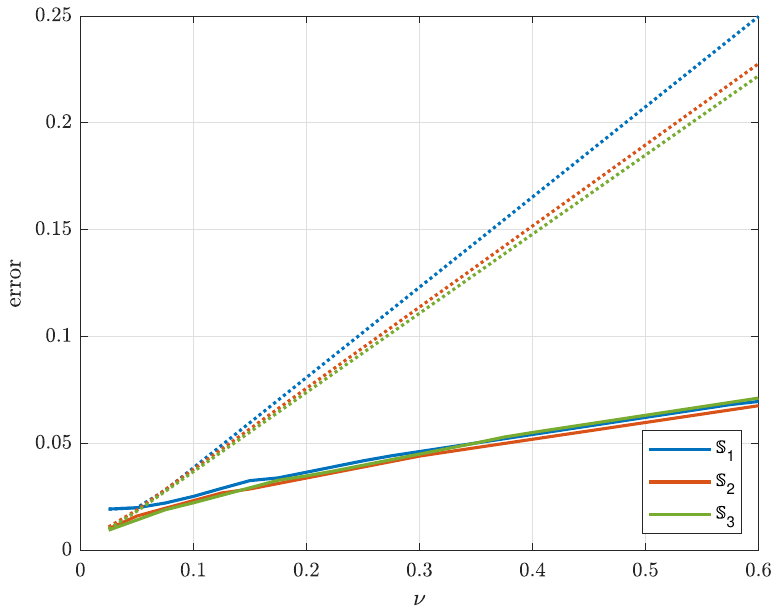}
\end{subfigure}

\caption{The approximation error with respect to the penalty parameter $\lambda$ and noise magnitude $\nu$ for the problem described in Example~\ref{ex:pdl2_franke}.}
\label{fig:pdl2_franke}
\end{figure}

\begin{example}
\label{ex:pdl2_franke}
As in Example \ref{ex:dl2_franke} we approximate the function \eqref{eq:franke_fun}, but this time we corrupt the function values computed on the $201 \times 201$ uniform grid of points by adding random noise drawn from the standard uniform distribution on the interval $(-\nu, \nu)$ for $\nu > 0$. We consider approximations obtained by the penalized method over the triangulation shown in Figure~\ref{fig:triangulation_sequence} (left). Figure~\ref{fig:pdl2_franke} (left) shows the approximation errors in dependence of $\lambda \in [10^{-8}, 10^2]$ with $\nu$ fixed to $0.25$. The dotted lines correspond to the approximations obtained for uncorrupted data ($\nu = 0$). On the other hand, Figure~\ref{fig:pdl2_franke} (right) shows the approximation errors in dependence of $\nu \in [0.025, 0.6]$ where for each $\nu$ and each spline space the value of $\lambda$ is computed in advance to be near-optimal. The dotted lines correspond to the approximations obtained without penalization ($\lambda = 0$). Note that on the unit square the function \eqref{eq:franke_fun} attains values between $0$ and $1.22$ and thus the values of $\nu$ in this test represent magnitudes of noise approximately from $2$ to $50$ percent of the function range.
\end{example}

\subsection{Approximations to boundary value problems}
\label{sec:boundary_value_problems}

Let $\Omega \subset \R^2$ be a connected bounded domain that can be parametrized by a $\C{1}$-smooth mapping $P \in \mathbb{S}_\ell \times \mathbb{S}_\ell$ defined on a polygonal domain $\Theta \subset \R^2$ triangulated by $\triangle$. We assume $P$ is a bijection with a nonvanishing Jacobian on $\Theta$ that maps the boundary $\partial \Theta$ of $\Theta$ to the boundary $\partial \Omega$ of $\Omega$.

\begin{figure}[!t]
\centering

\begin{subfigure}[b]{0.32\textwidth}
\centering
\includegraphics[width=\textwidth]{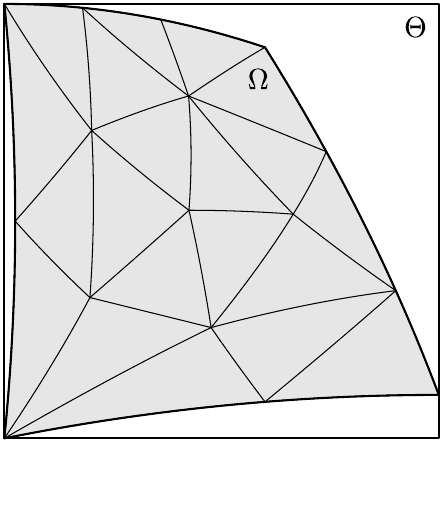}
\end{subfigure}
\hfill
\begin{subfigure}[b]{0.32\textwidth}
\centering
\includegraphics[width=\textwidth]{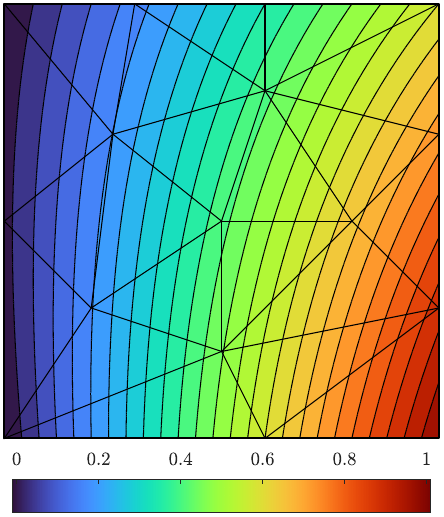}
\end{subfigure}
\hfill
\begin{subfigure}[b]{0.32\textwidth}
\centering
\includegraphics[width=\textwidth]{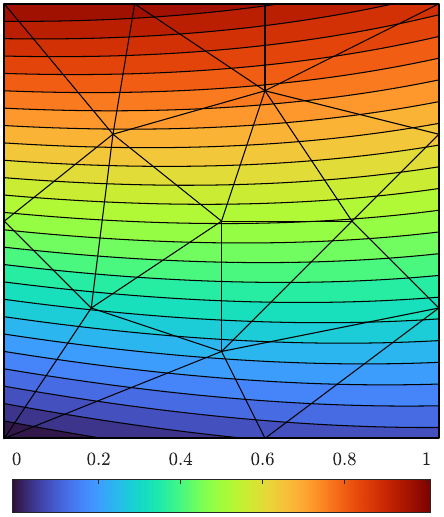}
\end{subfigure}

\caption{A domain with curved boundary $\Omega$ (left) discussed in Example~\ref{ex:geometric_mapping}, which is contained in the unit square $\Theta$. The domain $\Omega$ can be parametrized by a mapping whose components are quadratic polynomials (the contour plots of the first and the second component are shown in the middle and on the right, respectively). The left figure also depicts the partition of the domain that is obtained by using this mapping to map the triangulation from Example~\ref{ex:triangulation}.}
\label{fig:geometric_mapping}
\end{figure}

\begin{example}
\label{ex:geometric_mapping}
Figure~\ref{fig:geometric_mapping} (left) shows a domain $\Omega$ that can be exactly parametrized on the unit square $\Theta$ by a geometry mapping $P \in \poly_2 \times \poly_2$. Both components of this mapping (shown in the middle and on the right of Figure~\ref{fig:geometric_mapping}) are contained in $\mathbb{S}_\ell$, $\ell \in \set{1, 2, 3}$, and can be represented by the basis functions $S_1, S_2, \ldots, S_N$ with the help of \eqref{eq:bfunc} and \eqref{eq:brep}.
\end{example}

For given functions $F \in L^2(\Omega)$ and $G \in L^2(\partial \Omega)$ we consider the boundary value problem governed by the Poisson equation and the Dirichlet boundary conditions. The goal is to find the function $U$ satisfying
\begin{equation}
\label{eq:poisson_problem}
\left\{
\begin{aligned}
-\nabla^2 U(\bfm{p}) &= F(\bfm{p}) && \text{if} \ \bfm{p} \in \Omega \setminus \partial \Omega, \\
U(\bfm{p}) &= G(\bfm{p}) && \text{if} \ \bfm{p} \in \partial \Omega.
\end{aligned}
\right.
\end{equation}
In the isoparametric finite element method based on $\mathbb{S}_\ell$ we search for an approximation to $U$ of the form $S \circ P^{-1}$, $S \in \mathbb{S}_\ell$. After the Galerkin discretization in $\mathrm{span}(S_1 \circ P^{-1}, S_2 \circ P^{-1}, \ldots, S_n \circ P^{-1})$ this amounts to finding $S$ that admits
\begin{equation*}
\int_\Omega \nabla (S \circ P^{-1}) \bfm{\cdot} \nabla (S_j \circ P^{-1}) \, \mathrm{d}\Omega = \int_\Omega F \, (S_j \circ P^{-1}) \, \mathrm{d}\Omega, \quad j = 1, 2, \ldots, n,
\end{equation*}
or, equivalently,
\begin{equation}
\label{eq:isofem_conditions}
\int_\Theta \grad{S} \bfm{\cdot} \bfm{J}_P^{-1} \bfm{J}_P^{-T} \bfm{\cdot} \grad{S_j} \abs{\det(\bfm{J}_P)} \, \mathrm{d}\Theta = \int_\Theta (F \circ P) \, S_j \, \abs{\det(\bfm{J}_P)} \, \mathrm{d}\Theta, \quad j = 1, 2, \ldots, n,
\end{equation}
where $\bfm{J}_P$ is the Jacobian matrix of $P$.

We modify the expression \eqref{eq:srep} of $S$ to
\begin{equation*}
S = S_0 + \sum_{k = 1}^n c_k S_k, \quad S_0 = \sum_{k = n+1}^N c_k S_k.
\end{equation*}
First, we compute the coefficients $c_k$, $k = n+1, n+2, \ldots, N$, so that $S_0$ is the discrete best $L^2$ approximation to $G \circ P$ for a dense set of points on the boundary of $\Theta$. Then we transform \eqref{eq:isofem_conditions} into the system $\bfm{A} \bfm{\cdot} \bfm{c} = \bfm{b}$ determined by
\begin{equation*}
a_{j,k}^n = \int_\Theta \grad{S_k} \bfm{\cdot} \bfm{J}_P^{-1} \bfm{J}_P^{-T} \bfm{\cdot} \grad{S_j} \abs{\det(\bfm{J}_P)} \, \mathrm{d}\Theta, \quad
b_j^n = \int_\Theta \bracket{(F \circ P) \, S_j - \grad{S_0} \bfm{\cdot} \bfm{J}_P^{-1} \bfm{J}_P^{-T} \bfm{\cdot} \grad{S_j}}  \abs{\det(\bfm{J}_P)} \, \mathrm{d}\Theta.
\end{equation*}
The matrix $\bfm{A}$ is the stiffness matrix associated with the functions $S_1, S_2, \ldots, S_n$.

\begin{figure}[!t]
\centering

\begin{subfigure}[b]{0.49\textwidth}
\centering
\includegraphics[width=\textwidth]{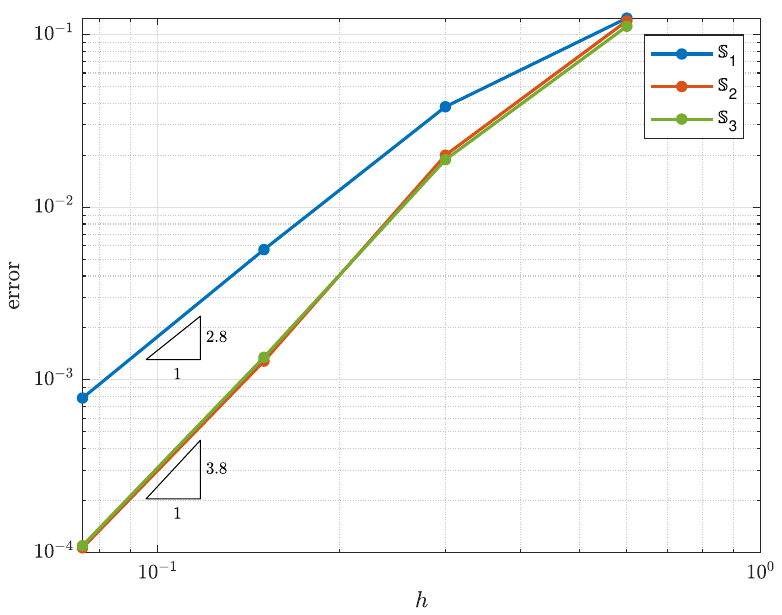}
\caption{Errors for Example~\ref{ex:fem}.}
\end{subfigure}
\hfill
\begin{subfigure}[b]{0.49\textwidth}
\centering
\includegraphics[width=\textwidth]{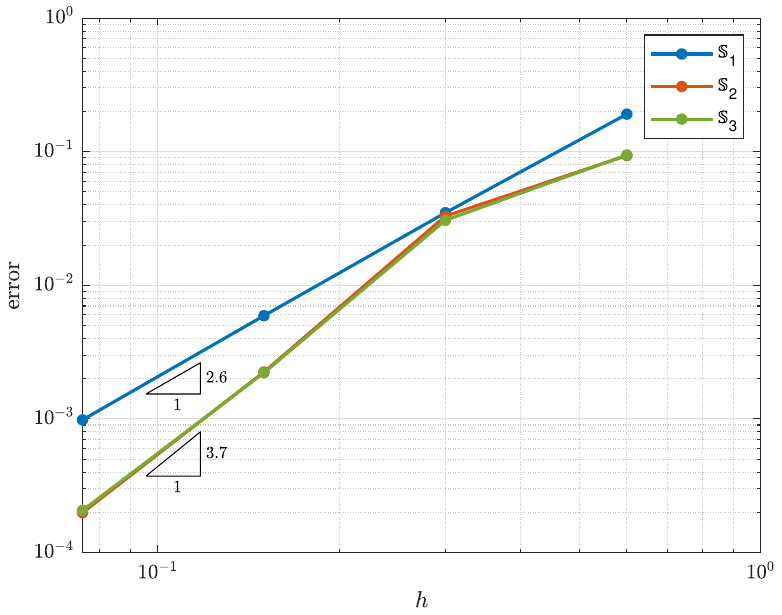}
\caption{Errors for Example~\ref{ex:isofem}.}
\end{subfigure}

\caption{Errors of approximations from $\mathbb{S}_\ell$, $\ell = 1, 2, 3$, for the problems discussed in Examples \ref{ex:fem} and \ref{ex:isofem} with respect to the length $h$ of the longest edge in the triangulation.}
\label{fig:fem_errors}
\end{figure}

\begin{example}
\label{ex:fem}
Let $P(x,y) = (x,y)$ and $\Omega = \Theta$ the unit square triangulated by the sequence of triangulations discussed in Example~\ref{ex:triangulation_sequence}. Applying the finite element method, we approximate the solution to the problem \eqref{eq:poisson_problem} with the manufactured solution
\begin{equation}
\label{eq:sin_fun}
U(x,y) = \sin\bracket{2 \pi (1-x) (1-y)}.
\end{equation}
The results are provided in Table~\ref{tab:fem} and in Figure~\ref{fig:fem_errors} (left).
\end{example}

\begin{table}[!t]
\centering
\begin{tabular}{|c|cc|cc|cc|}
\hline
& \multicolumn{2}{c|}{$\mathbb{S}_1$} & \multicolumn{2}{c|}{$\mathbb{S}_2$} & \multicolumn{2}{c|}{$\mathbb{S}_3$} \\
$n$ & error & cond & error & cond & error & cond \\ \hline 
24 & 1.25e-01 & 8.04e+00  & 1.20e-01 & 1.20e+01  & 1.12e-01 & 1.27e+01 \\
109 & 3.82e-02 & 3.82e+01  & 2.00e-02 & 4.35e+01  & 1.89e-02 & 4.55e+01 \\
459 & 5.69e-03 & 1.55e+02  & 1.28e-03 & 1.71e+02  & 1.35e-03 & 1.79e+02 \\
1879 & 7.85e-04 & 6.20e+02  & 1.06e-04 & 6.80e+02  & 1.09e-04 & 7.11e+02 \\
\hline
\end{tabular}
\caption{Approximation errors and condition numbers for the problem described in Example~\ref{ex:fem}.}
\label{tab:fem}
\end{table}

\begin{figure}[!t]
\centering

\begin{subfigure}[b]{0.24\textwidth}
\centering
\includegraphics[width=\textwidth]{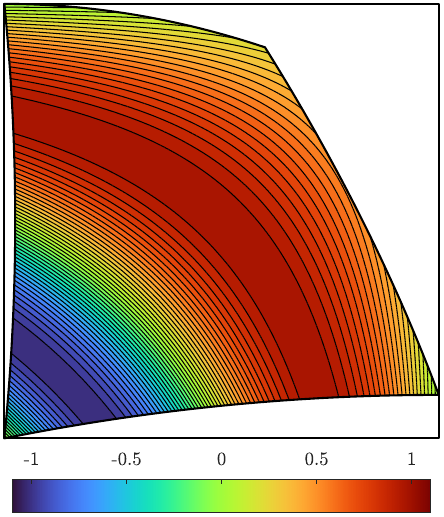}
\caption{Function to approximate.}
\end{subfigure}
\hfill
\begin{subfigure}[b]{0.24\textwidth}
\centering
\includegraphics[width=\textwidth]{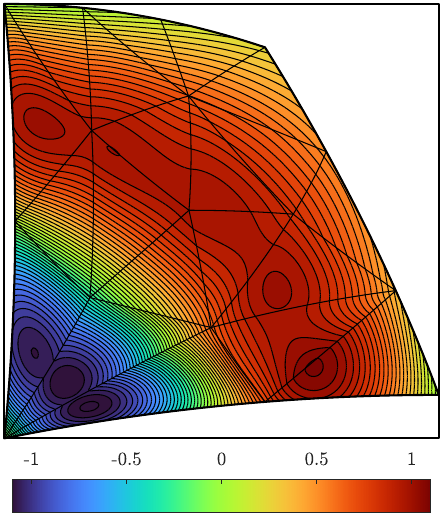}
\caption{Approximation from $\mathbb{S}_1$.}
\end{subfigure}
\hfill
\begin{subfigure}[b]{0.24\textwidth}
\centering
\includegraphics[width=\textwidth]{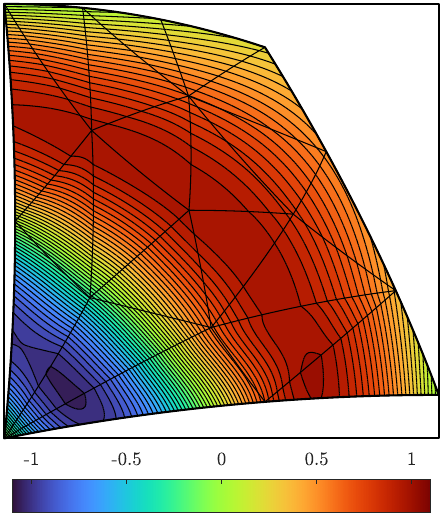}
\caption{Approximation from $\mathbb{S}_2$.}
\end{subfigure}
\hfill
\begin{subfigure}[b]{0.24\textwidth}
\centering
\includegraphics[width=\textwidth]{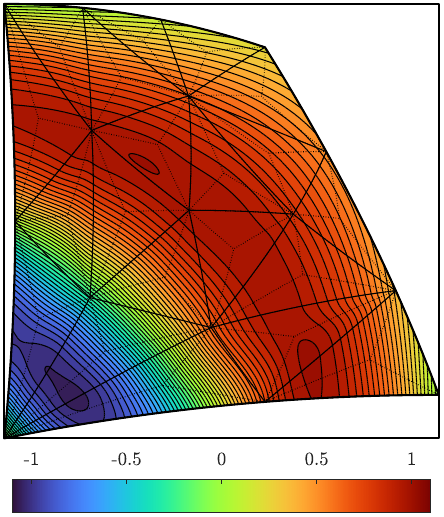}
\caption{Approximation from $\mathbb{S}_3$.}
\end{subfigure}

\caption{Contour plots of the function \eqref{eq:sin_fun} and its approximations from $\mathbb{S}_\ell$, $\ell = 1, 2, 3$, which are obtained by the isoparametric finite element method as discussed in Example~\ref{ex:isofem}.}
\label{fig:isofem}
\end{figure}

\begin{example}
\label{ex:isofem}
We repeat Example~\ref{ex:fem} for the mapping $P$ and domains $\Omega$ and $\Theta$ as defined in Example~\ref{ex:geometric_mapping}. The results are provided in Table~\ref{tab:isofem} and Figure~\ref{fig:fem_errors} (right). Here, the approximation errors are computed based on the $401 \times 401$ uniform grid of points mapped by $P$. Figure~\ref{fig:isofem} shows the contour plots of the approximations obtained for the initial triangulation.
\end{example}

\begin{table}[!t]
\centering
\begin{tabular}{|c|cc|cc|cc|}
\hline
& \multicolumn{2}{c|}{$\mathbb{S}_1$} & \multicolumn{2}{c|}{$\mathbb{S}_2$} & \multicolumn{2}{c|}{$\mathbb{S}_3$} \\
$n$ & error & cond & error & cond & error & cond \\ \hline 
24 & 1.91e-01 & 8.63e+00  & 9.36e-02 & 9.85e+00  & 9.45e-02 & 1.05e+01 \\
109 & 3.49e-02 & 3.36e+01  & 3.29e-02 & 3.54e+01  & 3.07e-02 & 3.74e+01 \\
459 & 5.92e-03 & 1.44e+02  & 2.24e-03 & 1.40e+02  & 2.22e-03 & 1.48e+02 \\
1879 & 9.83e-04 & 5.98e+02  & 1.99e-04 & 5.86e+02  & 2.07e-04 & 6.10e+02 \\
\hline
\end{tabular}
\caption{Approximation errors and condition numbers for the problem described in Example~\ref{ex:isofem}.}
\label{tab:isofem}
\end{table}

In \cite{ipbm_schumaker_19} an interesting alternative to the finite element method is proposed. As in the finite element method, the problem \eqref{eq:poisson_problem} is observed in the weak formulation but by considering all basis functions as test functions. The boundary conditions are imposed approximately with a penalty parameter, hence the name immersed penalty boundary method.

Instead of a parametric mapping between $\Theta$ and $\Omega$, we choose a polygonal domain $\Theta$ in such a way that it contains $\Omega$. Also, we choose $M \in \N$ points $\bfm{p}_d$, $d = 1, 2, \ldots, M$, on the boundary of $\Omega$. Moreover, we opt for a penalty parameter $\lambda > 0$ to ensure that $S$ approximately satisfies the boundary conditions and a penalty parameter $\mu > 0$ to ensure that $S$ is approximately $\C{1}$-smooth. Then we search for an approximation $S \in \mathbb{S}_\ell$ that solves
\begin{equation*}
\underset{S \in \mathbb{S}_\ell}{\mathrm{arg\,min}} \bracket{\sum_{j = 1}^N \bracket{\int_\Theta \bracket{F + \nabla^2 S} S_j \, \mathrm{d}\Theta}^2 + \lambda \sum_{d = 1}^M \bracket{G(\bfm{p}_d) - S(\bfm{p}_d)}^2 + \mu \sum_{l = 1}^{n_e - n_e^b} \xi_l(S)^2},
\end{equation*}
where $\xi_l$ is as defined in \eqref{eq:c1_condition_edge} and we impose (for notational convenience) that the edges in $E \setminus E^b$ are enumerated by the indices from $1$ to $n_e - n_e^b$. As $S \in \mathbb{S}_3$ is $\C{1}$-smooth, the third sum in the objective function should be ignored for $\ell = 3$. The solution to the minimization problem expressed in the form \eqref{eq:srep} can be obtained by solving the system $\bfm{A} \, \bfm{c} = \bfm{b}$ for
\begin{equation*}
\bfm{A} = \bfm{L}^T \bfm{L} + \lambda \, \bfm{S}^T \bfm{S} + \mu \bfm{E}^T \bfm{E}, \quad \bfm{b} = -\bfm{L}^T \bfm{f} + \lambda \, \bfm{S}^T \bfm{g},
\end{equation*}
where
\begin{itemize}
\item $\bfm{L} \in \R^{N \times N}$ is determined by the element $\int_\Theta (\nabla^2 S_k) \, S_j \, \mathrm{d}\Theta$ at position $(j,k)$,
\item $\bfm{S} \in \R^{M \times N}$ is determined by the element $S_k(\bfm{p}_d)$ at position $(d,k)$,
\item $\bfm{E} \in \R^{(n_e - n_e^b) \times N}$ is determined by the element $\xi_l(S_k)$ at position $(l,k)$,
\item $\bfm{f} \in \R^N$ is determined by the element $\int_\Theta F \, S_j \, \mathrm{d}\Theta$ at position $j$, and
\item $\bfm{g} \in \R^M$ is determined by the element $G(\bfm{p}_l)$ at position $l$.
\end{itemize}

\begin{figure}[!t]
\centering

\begin{subfigure}[b]{0.49\textwidth}
\centering
\includegraphics[width=\textwidth]{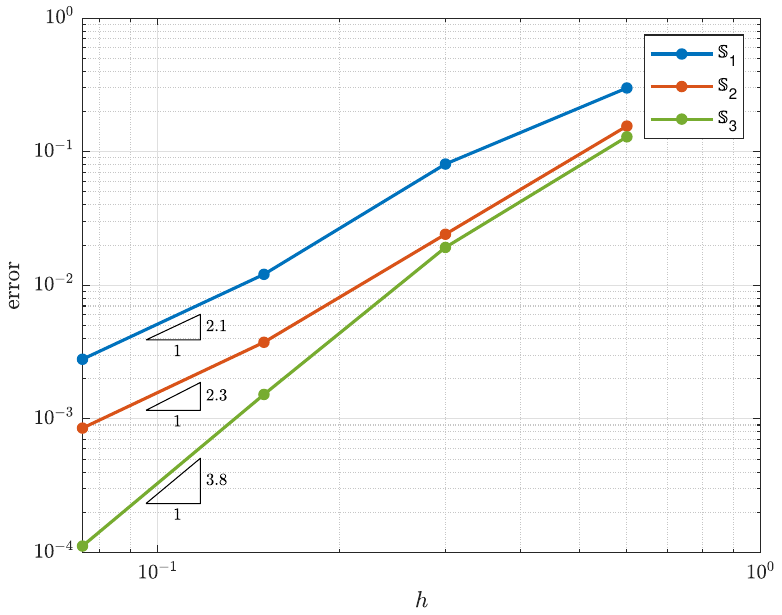}
\caption{Errors for Example~\ref{ex:ipbm}.}
\end{subfigure}
\hfill
\begin{subfigure}[b]{0.49\textwidth}
\centering
\includegraphics[width=\textwidth]{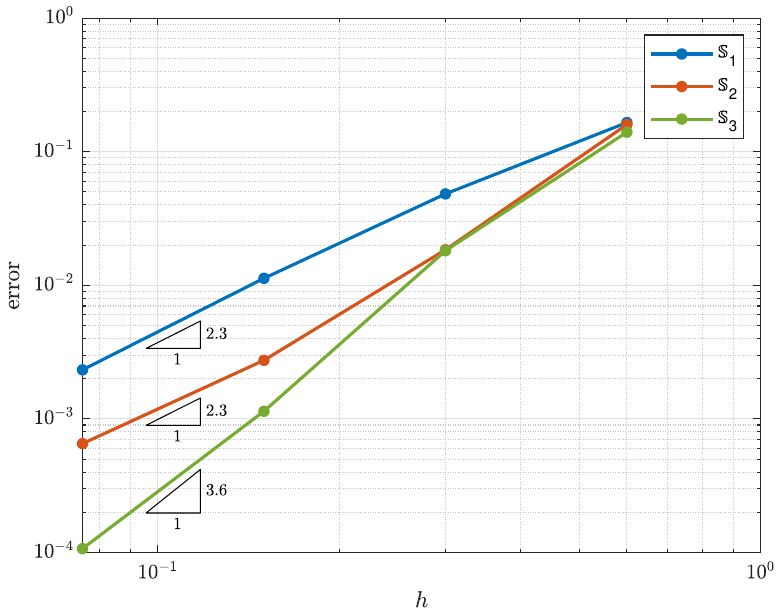}
\caption{Errors for Example~\ref{ex:ipbm_curved}.}
\end{subfigure}

\caption{Errors of approximations from $\mathbb{S}_\ell$, $\ell = 1, 2, 3$, for the problems discussed in Examples \ref{ex:ipbm} and \ref{ex:ipbm_curved} with respect to the length $h$ of the longest edge in the triangulation.}
\label{fig:ipbm_errors}
\end{figure}

\begin{example}
\label{ex:ipbm}
We consider the problem from Example~\ref{ex:fem} and solve it by using the immersed penalized boundary method on the unit square. We opt for $\lambda = 1$, $\mu = 1$, and $M = 800$ with the points $\bfm{p}_d$ spaced uniformly along the boundary of the domain. The results are provided in Table~\ref{tab:ipbm} and Figure~\ref{fig:ipbm_errors} (left).
\end{example}

\begin{table}[!t]
\centering
\begin{tabular}{|c|cc|cc|cc|}
\hline
& \multicolumn{2}{c|}{$\mathbb{S}_1$} & \multicolumn{2}{c|}{$\mathbb{S}_2$} & \multicolumn{2}{c|}{$\mathbb{S}_3$} \\
$N$ & error & cond & error & cond & error & cond \\ \hline 
48 & 3.00e-01 & 4.52e+02  & 1.56e-01 & 2.88e+02  & 1.29e-01 & 2.57e+02 \\
153 & 8.09e-02 & 4.44e+03  & 2.41e-02 & 4.33e+03  & 1.92e-02 & 3.51e+03 \\
543 & 1.21e-02 & 5.67e+04  & 3.75e-03 & 6.36e+04  & 1.52e-03 & 5.26e+04 \\
2043 & 2.80e-03 & 8.15e+05  & 8.53e-04 & 9.57e+05  & 1.12e-04 & 8.23e+05 \\
\hline
\end{tabular}
\caption{Approximation errors and condition numbers for the problem described in Example~\ref{ex:ipbm}.}
\label{tab:ipbm}
\end{table}

\begin{figure}[!t]
\centering

\begin{subfigure}[b]{0.24\textwidth}
\centering
\includegraphics[width=\textwidth]{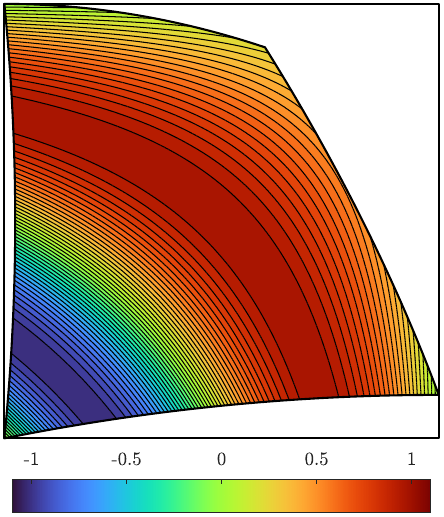}
\caption{Function to approximate.}
\end{subfigure}
\hfill
\begin{subfigure}[b]{0.24\textwidth}
\centering
\includegraphics[width=\textwidth]{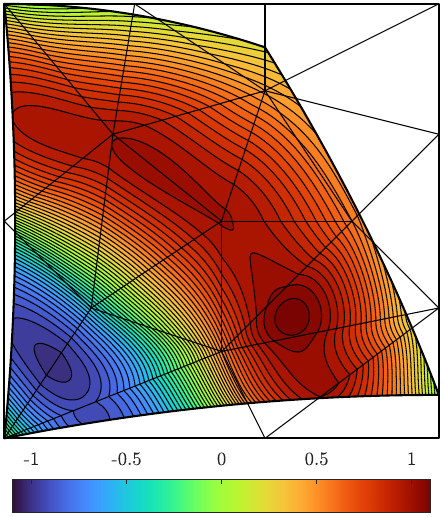}
\caption{Approximation from $\mathbb{S}_1$.}
\end{subfigure}
\hfill
\begin{subfigure}[b]{0.24\textwidth}
\centering
\includegraphics[width=\textwidth]{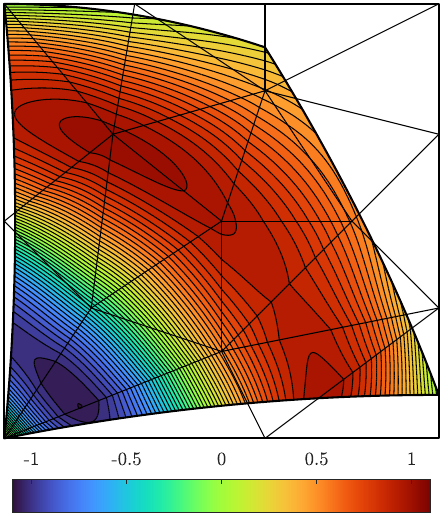}
\caption{Approximation from $\mathbb{S}_2$.}
\end{subfigure}
\hfill
\begin{subfigure}[b]{0.24\textwidth}
\centering
\includegraphics[width=\textwidth]{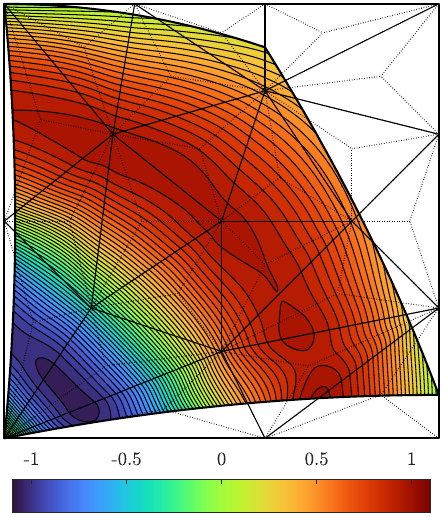}
\caption{Approximation from $\mathbb{S}_3$.}
\end{subfigure}

\caption{Contour plots of the function \eqref{eq:sin_fun} and its approximations from $\mathbb{S}_\ell$, $\ell = 1, 2, 3$, which are obtained by the immersed penalized boundary method as discussed in Example~\ref{ex:ipbm_curved}.}
\label{fig:ipbm_curved}
\end{figure}

\begin{example}
\label{ex:ipbm_curved}
This is a counterpart to Example~\ref{ex:isofem}. The domain $\Omega$ presented in Example~\ref{ex:geometric_mapping} is contained in the unit square $\Theta$. To apply the immersed penalized boundary method we use the same parameters $\lambda$, $\mu$, and $M$ as in Example~\ref{ex:ipbm} but map the uniformly spaced points from the boundary of the unit square to the boundary of $\Omega$ by $P$ to obtain $\bfm{p}_d \in \partial \Omega$. The results are provided in Table~\ref{tab:ipbm_curved} and Figure~\ref{fig:ipbm_errors} (right). The approximation errors are computed in the same way as in Example~\ref{ex:isofem}.
Figure~\ref{fig:ipbm_curved} shows contour plots of the approximations over the initial triangulation.
\end{example}

\begin{table}[!t]
\centering
\begin{tabular}{|c|cc|cc|cc|}
\hline
& \multicolumn{2}{c|}{$\mathbb{S}_1$} & \multicolumn{2}{c|}{$\mathbb{S}_2$} & \multicolumn{2}{c|}{$\mathbb{S}_3$} \\
$N$ & error & cond & error & cond & error & cond \\ \hline 
48 & 3.00e-01 & 4.52e+02  & 1.56e-01 & 2.88e+02  & 1.29e-01 & 2.57e+02 \\
153 & 8.09e-02 & 4.44e+03  & 2.41e-02 & 4.33e+03  & 1.92e-02 & 3.51e+03 \\
543 & 1.21e-02 & 5.67e+04  & 3.75e-03 & 6.36e+04  & 1.52e-03 & 5.26e+04 \\
2043 & 2.80e-03 & 8.15e+05  & 8.53e-04 & 9.57e+05  & 1.12e-04 & 8.23e+05 \\
\hline
\end{tabular}
\caption{Approximation errors and condition numbers for the problem described in Example~\ref{ex:ipbm_curved}.}
\label{tab:ipbm_curved}
\end{table}

\begin{figure}[!t]
\centering

\begin{subfigure}[b]{0.49\textwidth}
\centering
\includegraphics[width=\textwidth]{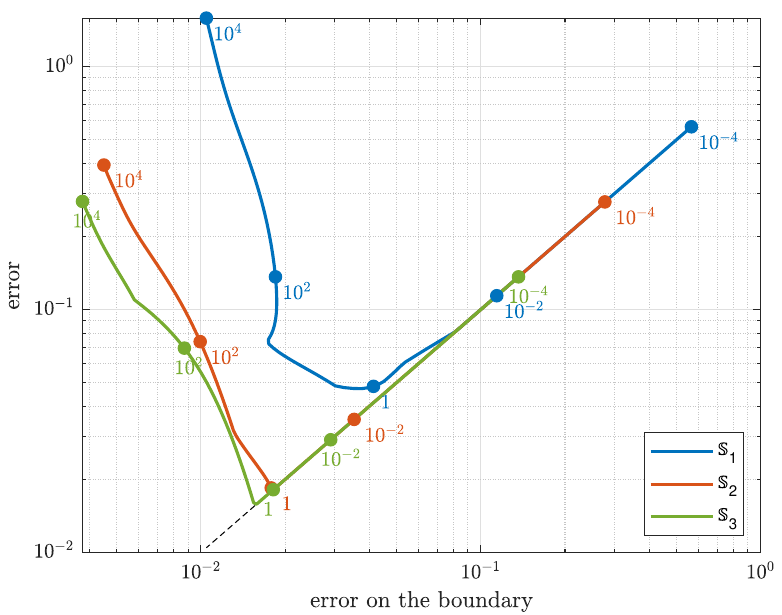}
\end{subfigure}
\hfill
\begin{subfigure}[b]{0.49\textwidth}
\centering
\includegraphics[width=\textwidth]{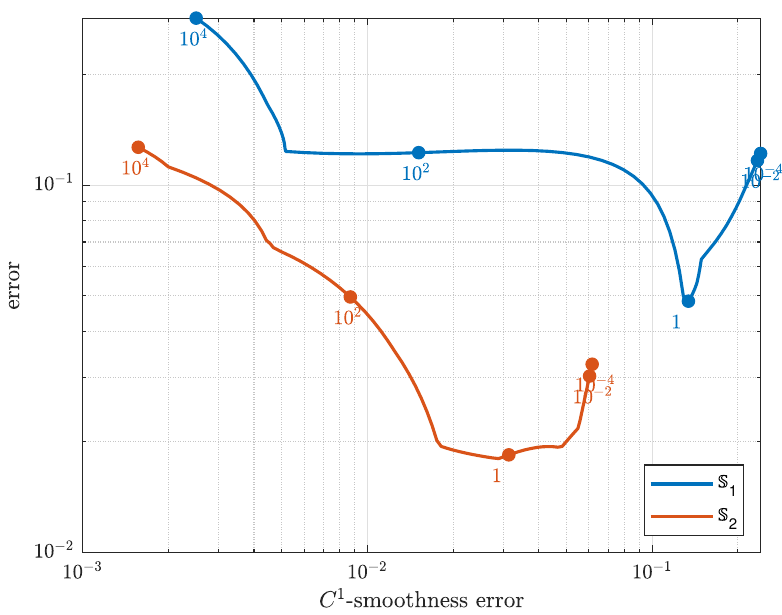}
\end{subfigure}

\caption{Impact of the penalization parameters $\lambda$ and $\mu$ on the accuracy of the approximations discussed in Example~\ref{ex:ipbm_curved_params}.
}
\label{fig:ipbm_curved_params}
\end{figure}

\begin{example}
\label{ex:ipbm_curved_params}
Finally, we report some more details on the choice of the penalization parameters $\lambda$ and $\mu$ in Example~\ref{ex:ipbm_curved}. We focus on the approximations over the triangulation shown in Figure~\ref{fig:triangulation_sequence} (left). Figure~\ref{fig:ipbm_curved_params} (left) justifies the choice $\lambda = 1$. For $\mu = 1$ it shows the relation between the approximation error on the entire domain and on the boundary when the parameter $\lambda$, which controls the importance of the boundary condition, takes the values between $10^{-4}$ and $10^4$. Figure~\ref{fig:ipbm_curved_params} (right) supports the choice $\mu = 1$. For $\lambda = 1$ it shows the relation between the approximation error and the $\C{1}$-smoothness error when the parameter $\mu$, which controls the importance of smoothness, takes the values between $10^{-4}$ and $10^4$. The $\C{1}$-smoothness error of an approximation $S$ is measured as the maximal absolute value of $\xi_l(S)$, $l = 1, 2, \ldots, n_e - n_e^b$. 
\end{example}

\section{Conclusion}
\label{sec:conclusion}

In the paper, three cubic macro-elements were considered that have the same characterization, \ie{}, they are determined by assigning three degrees of freedom to each vertex of the triangulation. However, the macro-elements differ in several aspects.

The first one is the classical Zienkiewicz element, which is on a triangle of the triangulation constructed by considering only the data corresponding to the three vertices of the triangle. It does not reproduce cubic polynomials and is only $\C{0}$-smooth, but it is simple to compute. The second one has cubic precision, which is achieved by taking advantage of the data corresponding to neighboring triangles. This considerably complicates the implementation but does not significantly increase the computational cost. However, the macro-element is still only $\C{0}$-smooth. The third approach guarantees both cubic precision and $\C{1}$-smoothness, thanks to using the underlying Clough--Tocher refinement of the triangulation. The disadvantage is that the refinement triples the number of triangles to work with and correspondingly slows the computations.

It was shown in the paper that despite different constructions, the considered macro-elements can be represented in a similar way by a basis consisting of three basis functions per vertex. The basis was utilized to numerically inspect the macro-elements in a variety of approximation methods. Tests indicate that in terms of stability the proposed representations behave very much alike. However, there is, as expected, no definite answer which of the three macro-elements performs best in terms of accuracy and computational cost, and the choice must be made with respect to the problem at hand. We finish with some remarks drawing the conclusions from the numerical experiments presented in Section~\ref{sec:examples}.

\begin{remark}
The lack of cubic precision implies nonoptimal convergence speed, which is a serious shortcoming when high-accuracy approximation is needed. This is clearly indicated in Examples \ref{ex:l2_franke}, \ref{ex:dl2_franke}, \ref{ex:fem}, and \ref{ex:isofem}.
\end{remark}

\begin{remark}
In the approximation of noisy data reported in Example~\ref{ex:pdl2_franke} all three macro-elements perform comparably. Rather surprisingly, neither the cubic precision nor the $\C{1}$-smoothness provide noteworthy advantage.
\end{remark}

\begin{remark}
For the considered problems of the best $L^2$ approximation (Examples \ref{ex:l2_franke} and \ref{ex:dl2_franke}) and solving the Poisson boundary problem via the finite element method (Examples \ref{ex:fem} and \ref{ex:isofem}), the two macro-elements with cubic precision provide practically identical accuracy. Nonetheless, the appearance of the approximations is considerably different as the $\C{0}$-smooth macro-element is prone to sharp transitions across the edges of the triangulation. This suggests the following. When only accuracy is important, it might be more efficient to use the $\C{0}$-smooth macro-element. If visual appearance matters, the results indicate that the $\C{1}$-smooth macro-element does not impair the quality of approximation.
\end{remark}

\begin{remark}
The lack of $\C{1}$-smoothness in the immersed penalized boundary method (Examples \ref{ex:ipbm} and \ref{ex:ipbm_curved}) clearly degrades the convergence rate of the approximations. The $\C{0}$-smooth approximations are less accurate than those obtained by the finite element method (Examples \ref{ex:fem} and \ref{ex:isofem}). The use of a penalty parameter to sanction the $\C{1}$-smoothness violation improves the visual appearance of the approximations (see Figure~\ref{fig:ipbm_curved}) but does not match the accuracy of the $\C{1}$-smooth approximations. In fact, the $\C{1}$-smooth macro-element performs remarkably well and provides as accurate results as in the finite element method.
\end{remark}

\section*{Acknowledgements}

The work of J.~Gro\v{s}elj was partially supported by the research programme P1-0294 of Javna agencija za znanstvenoraziskovalno in inovacijsko dejavnost Republike Slovenije (ARIS).

% \appendix

\bibliographystyle{elsarticle-num}
\bibliography{references}


%\begin{thebibliography}{00}



%\end{thebibliography}
\end{document}